\documentclass[11pt]{article}
\usepackage{amscd}
\usepackage{amssymb,amsfonts,amsmath,psfrag,amscd,cite}
\usepackage{amsmath}
  \usepackage{paralist}
  \usepackage{graphics}
  \usepackage{epsfig}
  \usepackage[colorlinks=true]{hyperref}
  \usepackage[all]{xy}
\usepackage{graphicx}
\newtheorem{theo}{Theorem}
\newtheorem{prop}[theo]{Proposition}

\newtheorem{defi}[theo]{Definition}

\newtheorem{rema}[theo]{Remark}

\makeatletter \@addtoreset{equation}{section}
\@addtoreset{theo}{section} \makeatother

\begin{document}
\date{}
\title{Symmetric Reduction of Regular Controlled Lagrangian\\
 System with Momentum Map}
\author{Hong Wang \\
School of Mathematical Sciences and LPMC,\\
Nankai University,  Tianjin 300071, P.R.China\\
E-mail: hongwang@nankai.edu.cn\\\\
March 12, 2021} \maketitle

{\bf Abstract.} In this paper, following the ideas in Marsden et
al.\cite{mawazh10}, we set up the regular reduction theory of a regular
controlled Lagrangian (RCL) system with symmetry and momentum map, by using
Legendre transformation and Euler-Lagrange vector field, and this
reduction is an extension of symmetric reduction theory of a regular
Lagrangian system under regular controlled Lagrangian equivalence
conditions. Considering the completeness of reduction, in order to
describe uniformly the RCL systems defined on a tangent bundle and on
its regular reduced spaces, we first define a kind of RCL systems on
a symplectic fiber bundle. Then we give a
good expression of the dynamical vector field of the RCL system,
such that we can describe the RCL-equivalence for the RCL systems.
Moreover, we introduce regular point and
regular orbit reducible RCL systems with symmetries and momentum
maps, by using the reduced Lagrange symplectic forms and the reduced Euler-Lagrange
vector fields, and prove the regular point
and regular orbit reduction theorems for the RCL systems and
regular Lagrangian systems, which explain
the relationships between RpCL-equivalence, RoCL-equivalence for
the reducible RCL systems with symmetries and RCL-equivalence for
the associated reduced RCL systems, as well as the relationship
of equivalences of the regular reducible Lagrangian systems, $R_p$-reduced
Lagrangian systems and $R_o$-reduced Lagrangian systems. \\



\tableofcontents

\section{Introduction}

It is well known that the theory of controlled mechanical systems
has formed an important subject in recent thirty years. Its research gathers
together some separate areas of research such as mechanics,
differential geometry and nonlinear control theory, etc., and the
emphasis of this research on geometry is motivated by the aim of
understanding the structure of equations of motion of the system in
a way that helps both for analysis and design. Thus, it is natural
to study the controlled mechanical systems by combining with the
analysis of dynamic systems and the geometric reduction theory of
Hamiltonian and Lagrangian systems.
Following the theoretical development of geometric
mechanics, a lot of important problems about this subject are being
explored and studied.\\

In 2005, we hope to study the mechanical system with control from
geometrical viewpoint. I and my students read the following two
papers of Professor Marsden and his students on a seminar of Nankai University,
see Chang et al. \cite{chbllemawo02, chma04}. We found that
there are some serious wrong of rigor for the definitions
of controlled Lagrangian (CL) system, controlled Hamiltonian (CH)
system and the reduced CH systems, the reduced CL system, as well as
CH-equivalence, CL-equivalence and the reduced CH-equivalence,
the reduced CL-equivalence in this two papers. In Marsden et al. \cite{mawazh10},
we have corrected and renewed carefully some of these wrong definitions
for the regular controlled Hamiltonian (RCH) system
and the reduced RCH systems. In this paper,
we shall consider the cases of the regular controlled Lagrangian
(RCL) system and the reduced RCL systems.\\

There are the following three aspects of these wrong
for CL system and the reduced CL systems in \cite{chbllemawo02,chma04}:\\

(1) The authors define CL system by using a wrong expression.
In fact, the CL system is defined in \cite{chbllemawo02,chma04},
by using the following expression
\begin{equation}
\mathcal{EL}(L)(q,\dot{q},\ddot{q})=F(q,\dot{q})+u(q,\dot{q}),
\label{1.1}
\end{equation}
where $\mathcal{EL}$ is the Euler-Lagrange operator and the bundle maps
$\mathcal{EL}(L): T^{(2)}Q \rightarrow T^*Q, $
and $F: TQ \rightarrow T^*Q,$ and the control
$u: TQ \rightarrow W (\subset T^*Q)$. This expression (1.1) can not be an equation,
because the left side of (1.1) is defined on the second order
tangent bundle $T^{(2)}Q$, and the right side of (1.1) is defined on the tangent
bundle $TQ$, and $T^{(2)}Q$ and $TQ$ are different spaces. Thus, it is
impossible to define the CL system by using a wrong expression. The
similar wrong appears in the definition 2.4 of \cite{chma04} for the
reduced CL system. In addition, it is worthy of noting that
the use of the above wrong expression (1.1)
has led to the wrong of the method of controlled Lagrangians to judge
the stabilization of mechanical systems,
see Bloch et al. \cite{bllema00, bllema01, blchlema01, blle02},
and this is also a serious problem should to be corrected carefully
in our next paper.\\

(2) The authors didn't consider the phase spaces of CL system and
the reduced CL system, that is, all of CL systems and the reduced
CL systems given in \cite{chbllemawo02,chma04}, have not the spaces
on which these systems are defined, see Definition 2.1 in
\cite{chbllemawo02} and Definition 2.1, 2.3 in \cite{chma04}. Thus,
it is impossible to give the actions of a Lie group
on the phase spaces of CL systems and their momentum maps,
also impossible to determine precisely the
reduced phase spaces of CL systems.\\

(3) The authors didn't consider the momentum map of the CL system with
symmetry, and didn't consider yet the change of geometrical
structures of the phase spaces of the CL systems,
and hence cannot determine precisely the geometrical
structures of phase spaces of the reduced CL systems. In fact,
it is not that all of CL systems in \cite{chma04} have same phase space $TQ$,
same action of Lie group $G$, and same reduced phase space $TQ /G$.
Different structures of geometry determine the different CL systems,
the the different reduced CL systems and their phase spaces.
Moreover, it is also impossible to give precisely the relationship of the
equivalences for the reduced CL systems, if don't consider
the different Lie group actions and momentum maps.\\

To sum up the above statement, we think that there
are a lot of serious wrong of rigor for the definitions
of CL system and its reduced CL systems, as well as
CL-equivalence and the reduced CL-equivalence
in Chang et al. \cite{chbllemawo02, chma04},
and we want to correct their work.
It is important to find these wrong, but the more important
is to correct well these wrong. It is worthy of noting that
we can not define directly a CL system on the second order tangent bundle $T^{(2)}Q$
or on the tangent bundle $TQ$ .
Because for $\mu \in \mathfrak{g}^\ast$, a regular value
of the momentum map $\mathbf{J}$, $G_\mu$ is the isotropy subgroup of the
coadjoint $G$-action at the point $\mu$,
the reduced second order tangent bundle $(T^{(2)}Q)/G_\mu$
or the reduced tangent bundle $(TQ)/G_\mu$ may not be a
second order tangent bundle or not be a tangent bundle of a configuration manifold.
If we define directly a CL system with symmetry on the second order tangent bundle $T^{(2)}Q$
or on the tangent bundle $TQ$,
then the reduced CL system may not have definition.
Thus, in order to set up the regular reduction theory for the CL system with symmetry,
in this paper, following the ideas in Marsden et al. \cite{mawazh10},
we have to correct and renew carefully these wrong definitions
in Chang et al. \cite{chbllemawo02, chma04}.\\

In this paper, we first consider that a regular controlled Lagrangian (RCL) system is a
regular Lagrangian system with external force and control.
In general, an RCL system under the action of external force and control is not
a regular Lagrangian system, however, it is a dynamical system closely related to a
regular Lagrangian system, and it can be explored and studied by extending
the methods for external force and control in the study of
regular Lagrangian systems. In consequence,
we can set up the regular reduction theory for an RCL system
with symmetry and momentum map,
by analyzing carefully the geometrical and topological
structures of the phase space and the reduced phase spaces of the
corresponding regular Lagrangian system.\\

A brief of outline of this paper is as follows. In the second
section, we review some relevant definitions and basic facts about
the regular Lagrangian system and its regular point and regular orbit reductions;
we also analyse the geometrical structures of the phase space and
the reduced phase spaces of the regular Lagrangian system,
which will be used in subsequent sections.
An RCL system is defined by using a (Lagrangian) symplectic form
on a symplectic fiber bundle and on the tangent bundle of a
configuration manifold, respectively,
and a good expression of the dynamical vector field for the RCL system
is given, and RCL-equivalence is
introduced in the third section. From the fourth section we begin to
discuss the RCL systems with symmetries and momentum maps
by combining with regular reduction theory of a regular Lagrangian system.
The regular point and regular orbit
reducible RCL systems are considered respectively in the fourth
section and the fifth section, and prove the regular point and
regular orbit reduction theorems for the RCL systems to explain the
relationships between the RpCL-equivalence, RoCL-equivalence for
the reducible RCL systems with symmetries and the RCL-equivalence for
the associated reduced RCL systems. We also study the relationship
of equivalences of the regular reducible Lagrangian systems, $R_p$-reduced
Lagrangian systems and $R_o$-reduced Lagrangian systems.
These research work develop the theory of symmetric reduction
for the RCL systems with symmetries and momentum maps,
and make us have much deeper
understanding and recognition for the structures of the regular
controlled mechanical systems.

\section{Legendre Transformation, Regular Lagrangian System
and Its Reduction}

In this section, we review some relevant definitions and basic facts
about Legendre transformation, the regular Lagrangian system and
its regular point and regular orbit reductions;
we also analyse the geometrical structures of the phase space and
the reduced phase spaces for a regular Lagrangian system,
which will be used in subsequent sections.
We shall follow the notations and conventions introduced in
Abraham and Marsden \cite{abma78},  Abraham et al.
\cite{abmara88}, Marsden \cite{ma92}, Marsden et al.
\cite{mamiorpera07}, Marsden and Ratiu \cite{mara99},
and Ortega and Ratiu \cite{orra04}.
For convenience, we assume that all manifolds in this paper are
real, smooth and finite dimensional. In particular, the following we always
assume that $Q$ is a smooth manifold with coordinates $q^i$, and
$TQ$ its tangent bundle with coordinates $(q^i,\dot{q}^i)$, and
$T^\ast Q$ its cotangent bundle with coordinates $(q^i,p_i)$, which
is the canonical cotangent bundle coordinates of $T^\ast Q$ and
$\theta_{0}=p_{i}\mathbf{d}q^{i}$ and $\omega_0= -\mathbf{d}\theta_{0}
=\mathbf{d}q^{i}\wedge \mathbf{d}p_{i}$ are canonical one-form and
canonical symplectic form on $T^{\ast}Q$, respectively,
where summation on repeated indices is understood.

\begin{defi}
Assume that $Q$ is an n-dimensional smooth manifold and the function
$L: TQ \rightarrow \mathbb{R}$. Then the map $\mathcal{F}L: TQ
\rightarrow T^{\ast}Q $ defined by
\begin{equation}
 \mathcal{F}L(v)w:=\frac{\mathrm{d}}{\mathrm{d}t}\mid_{t=0}L_{q}(v+tw),
 \;\;\; \forall  \; v ,w \in T_{q}Q ,
 \label{2.1}
\end{equation}
is fiber-preserving smooth map, which is called the fiber
derivative of $L$, where $L_{q}$ denotes the restrictions of $L$ to
the fiber over $q\in Q.$ If $\mathcal{F}L: TQ \rightarrow T^{\ast}Q
$ is a local diffeomorphism, then $L: TQ\rightarrow \mathbb{R}$ is
called a regular Lagrangian; and if $\mathcal{F}L:TQ \rightarrow
T^{\ast}Q$ is a diffeomorphism, then $L$ is called hyperregular.
\end{defi}

In the finite dimensional case, the local expression of the map
$\mathcal{F}L: TQ \rightarrow T^{\ast}Q $ is given by
\begin{equation}
 \mathcal{F}L(q^i, \dot{q}^i)=(q^i, \frac{\partial L}{\partial \dot{q}^i})
 =(q^i, p_i).
 \label{2.2}
\end{equation}
The change of data from $(q^i, \dot{q}^i)$ on $TQ$ to $(q^i,
p_i)$ on $T^{\ast}Q$, which is given by the map $\mathcal{F}L: TQ
\rightarrow T^{\ast}Q $, is called a {\bf Legendre transformation}. From
Marsden and Ratiu \cite{mara99}, we know that the Lagrangian $L$ is
regular, if the matrix $(\frac{\partial ^2 L}{\partial
\dot{q}^i\partial \dot{q}^j})$ is invertible. In the following by
using the Legendre transformation, we can give a definition of a
regular Lagrangian system as follows.

\begin{defi} (Regular Lagrangian System)
Assume that $Q$ is a smooth manifold, and $\theta_{0}$ and
$\omega_{0}$ are the canonical one form and the canonical symplectic
form on the cotangent bundle $T^{\ast}Q ,$ and the function $L: TQ
\rightarrow \mathbb{R}$ is hyperregular. Denote
$\theta^{L}:=(\mathcal{F}L)^{\ast}\theta_{0}$ and
$\omega^{L}:=(\mathcal{F}L)^{\ast}\omega_{0}$, where the bundle map
$(\mathcal{F}L)^{\ast}: T^*T^* Q \rightarrow T^*TQ. $ Then
$\theta^{L}$ and $\omega^{L}$ are called the Lagrangian one-form and
the Lagrangian symplectic form on the tangent bundle $TQ$,
respectively. Define an action $A: TQ \rightarrow \mathbb{R}$ given by
$A(v):=\mathcal{F}L(v)v, \; \forall v \in T_{q}Q $ and an energy $E_L: TQ
\rightarrow \mathbb{R}$ given by $E_L:= A-L$. If there exists a
vector field $\xi_{L}$ on $TQ$, such that the Euler-Lagrange
equation $\mathbf{i}_{\xi_{L}}\omega^{L}=\mathbf{d}E_L$ holds, then
$\xi_{L}$ is called an {\bf Euler-Lagrange vector field} of $L$, and the
triple $(TQ,\omega^L, L)$ is called a {\bf regular Lagrangian system}.
\end{defi}

In the finite dimensional case, the local expression of $\theta^{L}$
and $\omega^{L}$ are given by
$$\theta^{L}= \frac{\partial L}{\partial \dot{q}^i} dq^i,
\;\;\;\;\;\;\;\; \omega^{L}= \frac{\partial ^2 L}{\partial
\dot{q}^i\partial q^j}dq^i \wedge dq^j+ \frac{\partial ^2
L}{\partial \dot{q}^i\partial \dot{q}^j} dq^i \wedge d\dot{q}^j,$$
where summation on repeated indices is understood. Moreover, we know
that the energy $E_L$ is conserved along the flow of the
Euler-Lagrange vector field $\xi_L$, if $\xi_L$ satisfies a second order
equation, that is, $T\tau_Q \circ \xi_L = id_{TQ},$ where the map
$T\tau_Q: TTQ \rightarrow TQ,$ is the tangent map of the projection
$\tau_Q: TQ \rightarrow Q. $ Moreover, in a local coordinates of
$TQ$, an integral curve $(q(t), \dot{q}(t))$ of $\xi_L$ satisfies
the following Euler-Lagrange equations:
$$\frac{dq^i}{dt}= \dot{q}^i, \;\;\;\;\;\; \frac{d}{dt}\frac{\partial L}{\partial \dot{q}^i}
= \frac{\partial L}{\partial q^i}, \;\;\;\;\;\; i=1,2,\ldots,n.$$ If
$L$ is regular, then $\xi_L$ satisfies always the second order equation.\\

Furthermore, by using the Legendre transformation, the following
proposition gives a description of the equivalence between
the regular Lagrangian system $(TQ,\omega^L, L)$ and the Hamiltonian system
$(T^{\ast}Q, \omega_0, H)$ under the hyperregular case of $L$, see
Marsden and Ratiu \cite{mara99}.

\begin{prop}
Assume that $L:TQ\rightarrow \mathbb{R}$ is a hyperregular
Lagrangian on $TQ$, and define a function $H:=E_L\cdot
(\mathcal{F}L)^{-1}: T^{\ast}Q \rightarrow \mathbb{R}$. Then $H$ is
a hyperregular Hamiltonian on $T^{\ast}Q$, and the Hamiltonian
vector field $X_{H} \in TT^{\ast}Q$ and the Euler-Lagrange vector
field $\xi_{L} \in TTQ$ are $\mathcal{F}L$-related, i.e.
$T(\mathcal{F}L)\cdot \xi_{L}=X_{H}\cdot \mathcal{F}L, $ where
$T(\mathcal{F}L): TTQ \rightarrow TT^* Q$ is the tangent map of
$\mathcal{F}L: TQ \rightarrow T^{\ast}Q $, and the integral curves
of $\xi_{L}$ are mapped by $\mathcal{F}L$ onto integral curves of
$X_{H}$.
\end{prop}

It is well-known that Hamiltonian reduction theory is one of the
most active subjects in the study of modern analytical mechanics and
applied mathematics, in which a lot of deep and beautiful results
have been obtained, see the studies given by Abraham and Marsden
\cite{abma78},  Abraham et al. \cite{abmara88}, Arnold \cite{ar89},
Libermann and Marle \cite{lima87}, Marsden \cite{ma92}, Marsden et al.
\cite{mamiorpera07, mamora90}, Marsden and Perlmutter \cite{mape00},
Marsden and Ratiu \cite{mara99}, Marsden and
Weinstein \cite{mawe74}, Meyer \cite{me73},
Nijmeijer and Van der Schaft \cite {nivds90}and Ortega and Ratiu \cite{orra04},
in which the Marsden-Weinstein reduction
for the Hamiltonian system with symmetry and
momentum map is most important and foundational.
Now, for a regular Lagrangian system with symmetry and momentum map,
we can also give its regular point reduction as follows.\\

Let $Q$ be a smooth manifold and $TQ$ its tangent bundle with the
induced Lagrangian symplectic form $\omega^L$. Assume that $\Phi: G\times Q\rightarrow
Q$ be a smooth left action of a Lie group $G$ on $Q$, which is free
and proper, then the tangent lifted left action $\Phi^{T}: G\times
TQ\rightarrow TQ$ is also free, proper. Moreover, assume that
the action is symplectic with respect to
$\omega^L$, and admits an $\operatorname{Ad}^\ast$-equivariant
momentum map $\mathbf{J}_L: TQ\rightarrow \mathfrak{g}^\ast$, where
$\mathfrak{g}$ is the Lie algebra of $G$ and $\mathfrak{g}^\ast$ is
the dual of $\mathfrak{g}$. For a regular value of $\mathbf{J}_L$,
$\mu\in\mathfrak{g}^\ast$, denote $G_\mu=\{g\in
G|\operatorname{Ad}_g^\ast \mu=\mu \}$ the isotropy
subgroup of the co-adjoint $G$-action at the point
$\mu\in\mathfrak{g}^\ast$. Since $G_\mu (\subset G)$
acts freely and properly on $Q$ and on $TQ$, then $Q_\mu=Q/G_\mu$ is
a smooth manifold and the canonical projection
$\rho_\mu:Q\rightarrow Q_\mu$ is a surjective submersion. It follows
that $G_\mu$ acts also freely and properly on
$\mathbf{J}_L^{-1}(\mu)$, so that the space $(T
Q)_\mu=\mathbf{J}_L^{-1}(\mu)/G_\mu$ is a symplectic manifold with
the symplectic form $\omega^L_{\mu}$ uniquely characterized by the
relation
\begin{equation}\tau_\mu^\ast \cdot\omega^L_{\mu}=j_\mu^\ast
\cdot\omega^L. \label{2.3}\end{equation} The map
$j_\mu:\mathbf{J}_L^{-1}(\mu)\rightarrow TQ$ is the inclusion and
$\tau_\mu:\mathbf{J}_L^{-1}(\mu)\rightarrow (TQ)_\mu$ is the
projection. The pair $((TQ)_\mu,\omega^L_{\mu})$ is called the regular
point reduced space of $(TQ,\omega^L)$ at $\mu$.\\

Let $L: TQ\rightarrow \mathbb{R}$ be a $G$-invariant hyperregular Lagrangian, the
flow $F_t$ of the Euler-Lagrange vector field $\xi_L$ leaves the
connected components of $\mathbf{J}_L^{-1}(\mu)$ invariant and
commutes with the $G$-action, so it induces a flow $f_t^\mu$ on
$(TQ)_\mu$, defined by $f_t^\mu\cdot \tau_\mu=\tau_\mu \cdot
F_t\cdot j_\mu$, and the vector field $\xi_{l_\mu}$ generated by the
flow $f_t^\mu$ on $((TQ)_\mu,\omega^L_{\mu})$ is the reduced
Euler-Lagrange vector field with the associated regular point
reduced Lagrangian function $l_\mu: (TQ)_\mu \rightarrow \mathbb{R}$
defined by $l_\mu\cdot\tau_\mu=L\cdot j_\mu$, and the reduced
Euler-Lagrange equation
$\mathbf{i}_{\xi_{l_\mu}}\omega^L_{\mu}=\mathbf{d}E_{l_\mu}$ holds,
where the reduced energy $E_{l_\mu}: (TQ)_\mu \rightarrow
\mathbb{R}$ is given by $E_{l_\mu}:= A_\mu -l_\mu, $ and the reduced
action $A_\mu: (TQ)_\mu \rightarrow \mathbb{R} $ is given by
$A_\mu \cdot\tau_\mu= A\cdot j_\mu $,
and the Euler-Lagrange vector fields $\xi_L$ and $\xi_{l_\mu}$ are
$\tau_\mu$-related. Thus, we can introduce a kind of
regular point reducible Lagrangian systems as follows.
\begin{defi}
(Regular Point Reducible Lagrangian System) A 4-tuple $(TQ, G, \omega^L, L)$,
where the hyperregular Lagrangian $L: TQ \rightarrow \mathbb{R}$ is $G$-invariant,
is called a {\bf regular point reducible Lagrangian system}, if there exists a
point $\mu\in\mathfrak{g}^\ast$, which is a regular value of the
momentum map $\mathbf{J}_L$, such that the regular point reduced
system, that is, the 3-tuple $((TQ)_\mu, \omega^L_\mu,l_\mu)$,
where $(TQ)_\mu=\mathbf{J}^{-1}_L(\mu)/G_\mu$, $\tau_\mu^\ast
\cdot \omega^L_\mu= j_\mu^\ast \cdot \omega^L$, $l_\mu \cdot \tau_\mu= L \cdot
j_\mu$, is a regular Lagrangian system, which is simply written
as $R_p$-reduced Lagrangian system. Where $((TQ)_\mu,\omega^L_\mu)$ is the
$R_p$-reduced space, the function $l_\mu: (TQ)_\mu \rightarrow
\mathbb{R}$ is called the $R_p$-reduced Lagrangian.
\end{defi}

We know that the orbit reduction of a Hamiltonian system
is an alternative approach to symplectic reduction given
by Kazhdan, Kostant and Sternberg \cite{kakost78} and Marle \cite{ma76},
which is different from the Marsden-Weinstein reduction.
Now, for a regular Lagrangian system with symmetry and momentum map,
we can also give its regular orbit reduction as follows,
which is different from the above regular point reduction.\\

Assume that $\Phi:G\times Q\rightarrow Q$ is a smooth left action of a Lie
group $G$ on $Q$, which is free and proper, then the tangent lifted
left action $\Phi^{T}:G \times TQ\rightarrow TQ$ is also free and proper.
Moreover, assume that the action is
symplectic with respect to $\omega_L$, and admits an
$\operatorname{Ad}^\ast$-equivariant momentum map $\mathbf{J}_L:
TQ\rightarrow \mathfrak{g}^\ast$. For a regular value of the
momentum map $\mathbf{J}_L$, $\mu\in \mathfrak{g}^\ast$,
$\mathcal{O}_\mu=G\cdot \mu\subset \mathfrak{g}^\ast$
is the $G$-orbit of the co-adjoint $G$-action
through the point $\mu$. Since $G$ acts freely, properly and
symplectically on $TQ$ with respect to $\omega_L$, then the quotient space $(T
Q)_{\mathcal{O}_\mu}= \mathbf{J}_L^{-1}(\mathcal{O}_\mu)/G$ is a
regular quotient symplectic manifold with the reduced symplectic
form $\omega^L_{\mathcal{O}_\mu}$ uniquely characterized by the
relation
\begin{equation}
j_{\mathcal{O}_\mu}^\ast\cdot \omega^L=\tau_{\mathcal{O}_{\mu}}^\ast
\cdot\omega^L_{\mathcal{O}
_\mu}+(\mathbf{J}_{L})_{\mathcal{O}_\mu}^\ast\cdot\omega^{L+}_{\mathcal{O}_\mu},
\label{2.4}\end{equation}
where $(\mathbf{J}_L)_{\mathcal{O}_\mu}$
is the restriction of the momentum map $\mathbf{J}_L$ to
$\mathbf{J}_L^{-1}(\mathcal{O}_\mu)$, that is,
$(\mathbf{J}_L)_{\mathcal{O}_\mu}=\mathbf{J}_L\cdot
j_{\mathcal{O}_\mu}$ and $\omega_{\mathcal{O}_\mu}^{L+}$ and $\omega_{\mathcal{O}_\mu}^+$ are the
$+$-symplectic structures on the orbit $\mathcal{O}_\mu$ given by
\begin{equation}
\omega_{\mathcal{O}_\mu}^{L+}(\nu)(\xi, \; \eta)
=\omega_{\mathcal{O}_\mu}^
+(\nu)(\xi_{\mathfrak{g}^\ast}(\nu),\eta_{\mathfrak{g}^\ast}(\nu))
=<\nu,[\xi,\eta]>,\;\; \forall\;\nu\in\mathcal{O}_\mu, \;
\xi,\eta\in \mathfrak{g}, \; \xi_{\mathfrak{g}^\ast}, \eta_{\mathfrak{g}^\ast}
\in \mathfrak{g}^*. \label{2.5}\end{equation} The maps
$j_{\mathcal{O}_\mu}:\mathbf{J}_L^{-1}(\mathcal{O}_\mu)\rightarrow T
Q$ and
$\tau_{\mathcal{O}_\mu}:\mathbf{J}_L^{-1}(\mathcal{O}_\mu)\rightarrow
(TQ)_{\mathcal{O}_\mu}$ are natural injection and the projection,
respectively. The pair $((T
Q)_{\mathcal{O}_\mu},\omega^L_{\mathcal{O}_\mu})$ is called the
regular orbit reduced space of $(TQ,\omega^L)$ at the point $\mu$.\\

Let $L: TQ\rightarrow \mathbb{R}$ be a $G$-invariant hyperregular Lagrangian, the
flow $F_t$ of the Euler-Lagrange vector field $\xi_L$ leaves the
connected components of $\mathbf{J}_L^{-1}(\mathcal{O}_\mu)$
invariant and commutes with the $G$-action, so it induces a flow
$f_t^{\mathcal{O}_\mu}$ on $(TQ)_{\mathcal{O}_\mu}$, defined by
$f_t^{\mathcal{O}_\mu}\cdot
\tau_{\mathcal{O}_\mu}=\tau_{\mathcal{O}_\mu} \cdot F_t\cdot
j_{\mathcal{O}_\mu}$, and the vector field
$\xi_{l_{\mathcal{O}_\mu}}$ generated by the flow
$f_t^{\mathcal{O}_\mu}$ on $((T
Q)_{\mathcal{O}_\mu},\omega^L_{\mathcal{O}_\mu})$ is the reduced Euler-Lagrange
vector field with the associated regular orbit reduced Lagrangian
function $l_{\mathcal{O}_\mu}:(TQ)_{\mathcal{O}_\mu}\rightarrow
\mathbb{R}$ defined by $l_{\mathcal{O}_\mu}\cdot
\tau_{\mathcal{O}_\mu}= L\cdot j_{\mathcal{O}_\mu}$, and the reduced
Euler-Lagrange equation
$\mathbf{i}_{\xi_{l_{\mathcal{O}_\mu}}}\omega^L_{\mathcal{O}_\mu}=\mathbf{d}E_{l_{\mathcal{O}_\mu}}$
holds, where the reduced energy $E_{l_{\mathcal{O}_\mu}}:
(TQ)_{\mathcal{O}_\mu} \rightarrow \mathbb{R}$ given by
$E_{l_{\mathcal{O}_\mu}}:= A_{\mathcal{O}_\mu} -l_{\mathcal{O}_\mu},
$ and the reduced action $A_{\mathcal{O}_\mu}: (TQ)_{\mathcal{O}_\mu} \rightarrow
\mathbb{R}, $ given by $A_{\mathcal{O}_\mu}\cdot
\tau_{\mathcal{O}_\mu}= A\cdot j_{\mathcal{O}_\mu}$, and the Euler-Lagrange vector fields
$\xi_L$ and $\xi_{l_{\mathcal{O}_\mu}}$ are $\tau_{\mathcal{O}_\mu}$-related.
Thus, we can introduce a kind of the
regular orbit reducible Lagrangian systems as follows.
\begin{defi}
(Regular Orbit Reducible Lagrangian System) A 4-tuple $(TQ, G,
\omega^L,L)$, where the hyperregular Lagrangian
$L: TQ\rightarrow \mathbb{R}$ is $G$-invariant,
is called a {\bf regular orbit reducible Lagrangian
system}, if there exists an orbit $\mathcal{O}_\mu, \;
\mu\in\mathfrak{g}^\ast$, where $\mu$ is a regular value of the
momentum map $\mathbf{J}_L$, such that the regular orbit reduced
system, that is, the 3-tuple
$((TQ)_{\mathcal{O}_\mu},\omega^L_{\mathcal{O}_\mu},l_{\mathcal{O}_\mu})$, where
$(TQ)_{\mathcal{O}_\mu}=\mathbf{J}_L^{-1}(\mathcal{O}_\mu)/G$,
$\tau_{\mathcal{O}_\mu}^\ast \cdot \omega^L_{\mathcal{O}_\mu}
=j_{\mathcal{O}_\mu}^\ast \cdot \omega^L-(\mathbf{J}_L)_{\mathcal{O}_\mu}^\ast
\cdot \omega^{L+}_{\mathcal{O}_\mu}$, $l_{\mathcal{O}_\mu}\cdot
\tau_{\mathcal{O}_\mu} =L\cdot j_{\mathcal{O}_\mu}$,
is a regular Lagrangian system, which is
simply written as $R_o$-reduced Lagrangian system. Where
$((TQ)_{\mathcal{O}_\mu},\omega^L_{\mathcal{O}_\mu})$ is the
$R_o$-reduced space, the function
$l_{\mathcal{O}_\mu}:(TQ)_{\mathcal{O}_\mu}\rightarrow \mathbb{R}$
is called the $R_o$-reduced Lagrangian.
\end{defi}

In the following we shall give a precise analysis for
the geometrical structures of the regular point reduced space
$((TQ)_\mu,\omega^L_{\mu})$ and the regular orbit reduced space
$((TQ)_{\mathcal{O}_\mu},\omega^L_{\mathcal{O}_\mu})$.
Assume that the Lagrangian $L: TQ\rightarrow \mathbb{R}$ is
hyperregular, then the Legendre transformation $\mathcal{F}L: TQ
\rightarrow T^{\ast}Q $ is a diffeomorphism.
If the cotangent lift $G$-action $\Phi^{T*}: G\times T^\ast Q\rightarrow T^\ast Q$
is free, proper and symplectic with respect to the
canonical symplectic form $\omega_0$ on $T^*Q$, and has an
$\operatorname{Ad}^\ast$-equivariant momentum map $\mathbf{J}:T^\ast
Q\to \mathfrak{g}^\ast$ given by
$<\mathbf{J}(\alpha_q),\xi>=\alpha_q(\xi_Q(q)), $
where $\alpha_q \in T^*_q Q$ and $\xi\in
\mathfrak{g}$, $\xi_Q(q)$ is the value of the infinitesimal
generator $\xi_Q$ of the $G$-action at $q\in Q$, $<,>:
\mathfrak{g}^\ast \times \mathfrak{g}\rightarrow \mathbb{R}$ is the
duality pairing between the dual $\mathfrak{g}^\ast $ and
$\mathfrak{g}$. Then we have that the following theorem holds.
\begin{theo}
The momentum map $\mathbf{J}_L: TQ\rightarrow \mathfrak{g}^\ast$
given by $\mathbf{J}_L= \mathbf{J}\cdot \mathcal{F}L, $ is
$\operatorname{Ad}^\ast$-equivariant,
if the Lagrangian $L: TQ\rightarrow \mathbb{R}$ is
hyperregular, and the Legendre transformation $\mathcal{F}L: TQ
\rightarrow T^{\ast}Q $ is $(\Phi^T, \; \Phi^{T*})$-equivariant.
Moreover, if $\mu\in \mathfrak{g}^\ast $ is a regular
value of the momentum map $\mathbf{J}$, then
$\mu $ is also a regular
value of the momentum map $\mathbf{J}_L$.
\end{theo}

\noindent {\bf Proof:} We first prove that
the momentum map $\mathbf{J}_L: TQ\rightarrow \mathfrak{g}^\ast$
is $\operatorname{Ad}^\ast$-equivariant.
Since the Lagrangian $L: TQ\rightarrow \mathbb{R}$ is
hyperregular, then the Legendre transformation $\mathcal{F}L: TQ
\rightarrow T^{\ast}Q $ is a diffeomorphism. Because the momentum map
$\mathbf{J}: T^\ast Q\to \mathfrak{g}^\ast$ is
$\operatorname{Ad}^\ast$-equivariant, we have that
$\operatorname{Ad}^\ast \cdot \mathbf{J}= \mathbf{J}\cdot \Phi^{T*}$.
Note that the Legendre transformation $\mathcal{F}L: TQ
\rightarrow T^{\ast}Q $ is $(\Phi^T, \; \Phi^{T*})$-equivariant,
then we have that $\Phi^{T*}\cdot \mathcal{F}L= \mathcal{F}L \cdot \Phi^{T}$.
From the following commutative Diagram-1,
\[
\begin{CD}
TQ @> \Phi^T >> TQ @> \mathbf{J}_L >> \mathfrak{g}^\ast \\
@V \mathcal{F}L VV @V \mathcal{F}L VV @V \operatorname{Ad}^\ast VV \\
T^*Q @> \Phi^{T*} >> T^*Q @> \mathbf{J} >> \mathfrak{g}^\ast
\end{CD}
\]
$$\mbox{Diagram-1}$$
we can obtain that
$$
\operatorname{Ad}^\ast \cdot \mathbf{J}_L
= \operatorname{Ad}^\ast \cdot \mathbf{J}\cdot \mathcal{F}L
= \mathbf{J}\cdot \Phi^{T*} \cdot \mathcal{F}L
= \mathbf{J}\cdot \mathcal{F}L \cdot \Phi^{T}
= \mathbf{J}_L \cdot \Phi^{T}.
$$
Thus, the momentum map $\mathbf{J}_L: TQ\rightarrow \mathfrak{g}^\ast$
is $\operatorname{Ad}^\ast$-equivariant.\\

Next, if $\mu\in \mathfrak{g}^\ast $ is a regular
value of the momentum map $\mathbf{J}$, then there exists
an $\alpha \in T^*Q,$ such that $\mathbf{J}(\alpha)=\mu. $
Since the Legendre transformation $\mathcal{F}L: TQ
\rightarrow T^{\ast}Q $ is a diffeomorphism, we have that
$v= \mathcal{F}L^{-1}(\alpha)\in TQ,$ such that
$$
\mathbf{J}_L(v)= \mathbf{J}\cdot \mathcal{F}L( \mathcal{F}L^{-1}(\alpha))
= \mathbf{J}(\alpha)=\mu.$$
Thus, $\mu\in \mathfrak{g}^\ast $ is also a regular
value of the momentum map $\mathbf{J}_L$.
\hskip 1cm $\blacksquare$ \\

For a given $\mu\in\mathfrak{g}^\ast$, a regular value of
the momentum map $\mathbf{J}:T^\ast
Q\to \mathfrak{g}^\ast$, denote by $G_\mu$ the isotropy subgroup of the
co-adjoint $G$-action at the point $\mu$, then
the Marsden-Weinstein reduced space $(T^\ast
Q)_\mu=\mathbf{J}^{-1}(\mu)/G_\mu$ is a symplectic manifold with the
symplectic form $\omega_\mu$ uniquely characterized by the relation
\begin{equation}\pi_\mu^\ast \cdot \omega_\mu=i_\mu^\ast
\cdot \omega_0. \label{2.6}
\end{equation}
The map
$i_\mu:\mathbf{J}^{-1}(\mu)\rightarrow T^\ast Q$ is the inclusion
and $\pi_\mu:\mathbf{J}^{-1}(\mu)\rightarrow (T^\ast Q)_\mu$ is the
projection. From Marsden et al.\cite{mamiorpera07},
we know that the classification of symplectic
reduced spaces of a cotangent bundle is given as follows. (1) If $\mu=0$, the
symplectic reduced space of cotangent bundle $T^\ast Q$ at $\mu=0$
is given by $((T^\ast Q)_\mu, \omega_\mu)= (T^\ast(Q/G), \hat{\omega}_0)$,
where $\hat{\omega}_0$ is the canonical symplectic form of cotangent
bundle $T^\ast (Q/G)$. Thus, the symplectic reduced space $((T^\ast
Q)_\mu, \omega_\mu)$ at $\mu=0$ is a symplectic vector bundle. (2)
If $\mu\neq0$, and $G$ is Abelian, then $G_\mu=G$, in this case the
regular point symplectic reduced space $((T^*Q)_\mu, \omega_\mu)$ is
symplectically diffeomorphic to symplectic vector bundle $(T^\ast
(Q/G), \hat{\omega}_0-B_\mu)$, where $B_\mu$ is a magnetic term. (3) If
$\mu\neq0$, and $G$ is not Abelian and $G_\mu\neq G$, in this case
the regular point symplectic reduced space $((T^*Q)_\mu,
\omega_\mu)$ is symplectically diffeomorphic to a symplectic fiber
bundle over $T^\ast (Q/G_\mu)$ with fiber to be the co-adjoint orbit
$\mathcal{O}_\mu$, see the cotangent bundle reduction
theorem---bundle version, also see Marsden and Perlmutter
\cite{mape00}. Comparing the regular point reduced spaces
$((TQ)_\mu, \omega^L_\mu)$ and $((T^*Q)_\mu, \omega_\mu)$ at the point $\mu$,
we have that the following theorem holds.
\begin{theo}
Assume that the Lagrangian $L: TQ\rightarrow \mathbb{R}$ is
hyperregular, and the Legendre transformation $\mathcal{F}L: TQ
\rightarrow T^{\ast}Q $ is $(\Phi^T,\; \Phi^{T*})$-equivariant,
then the regular point reduced space $((TQ)_\mu, \omega^L_\mu)$
of $(TQ, \omega^L)$ at $\mu$ is symplectically diffeomorphic
to the regular point reduced space $((T^*Q)_\mu, \omega_\mu)$ of $(T^*Q, \omega_0)$ at $\mu$, and
hence is also symplectically diffeomorphic to a symplectic fiber
bundle.
\end{theo}

\noindent {\bf Proof:} Since the Lagrangian $L: TQ\rightarrow \mathbb{R}$ is
hyperregular, then the Legendre transformation $\mathcal{F}L: TQ
\rightarrow T^{\ast}Q $ is a diffeomorphism. Because $\mathcal{F}L$
is $(\Phi^T,\; \Phi^{T*})$-equivariant, that is,
$\Phi^{T*}\cdot \mathcal{F}L= \mathcal{F}L \cdot \Phi^{T}$,
then we can define a map
$(\mathcal{F}L)_\mu: (TQ)_\mu\rightarrow (T^*Q)_\mu $ given by
$(\mathcal{F}L)_\mu \cdot \tau_\mu= \pi_\mu \cdot \mathcal{F}L, $
and $i_\mu \cdot \mathcal{F}L= \mathcal{F}L \cdot j_\mu, $
see the following commutative Diagram-2,
which is well-defined and a diffeomorphism.
$$
\begin{CD}
\mathbf{J}_L^{-1}(\mu)\subset TQ @> \mathcal{F}L  >> \mathbf{J}^{-1}(\mu) \subset T^*Q \\
@V \tau_\mu VV @VV \pi_\mu  V \\
(TQ)_\mu @> (\mathcal{F}L)_\mu >> (T^*Q)_\mu
\end{CD}
$$
$$\mbox{Diagram-2}$$

We shall prove that
$(\mathcal{F}L)_\mu$ is symplectic, that is,
$(\mathcal{F}L)_\mu^* \cdot \omega_\mu= \omega^L_\mu .$
In fact, from (2.6) and (2.3), we have that
\begin{align*}
\tau_{\mu}^\ast \cdot (\mathcal{F}L)_\mu^* \cdot \omega_\mu
& =((\mathcal{F}L)_\mu \cdot \tau_{\mu})^*\cdot \omega_\mu
=( \pi_\mu \cdot \mathcal{F}L)^*\cdot \omega_\mu
= (\mathcal{F}L)^*\cdot  \pi_{\mu}^* \cdot \omega_\mu \\
& =(\mathcal{F}L)^*\cdot  i_{\mu}^* \cdot \omega_0
=(i_\mu \cdot \mathcal{F}L)^* \cdot \omega_0
=(\mathcal{F}L \cdot j_\mu)^*\cdot \omega_0\\
& =j_{\mu}^* \cdot (\mathcal{F}L)^*\cdot \omega_0
 =j_{\mu}^* \cdot \omega^L= \tau_{\mu}^\ast \cdot \omega^L_\mu.
\end{align*}
Notice that $\tau_{\mu}$ is surjective, and hence
$(\mathcal{F}L)_\mu^* \cdot \omega_\mu= \omega^L_\mu .$
Thus, the regular point reduced space $((TQ)_\mu, \omega^L_\mu)$
of $(TQ, \omega^L)$ at $\mu$
is symplectically diffeomorphic to the regular point reduced space $((T^*Q)_\mu, \omega_\mu)$
of $(T^*Q, \omega_0)$ at $\mu$.
From Marsden et al.\cite{mamiorpera07},
we know that the space $((T^*Q)_\mu, \omega_\mu)$
is symplectically diffeomorphic to a symplectic fiber
bundle, and hence $((TQ)_\mu, \omega^L_\mu)$
is also symplectically diffeomorphic to a symplectic fiber bundle.
\hskip 1cm $\blacksquare$ \\

For a given $\mu\in\mathfrak{g}^\ast$, a regular value of
the momentum map $\mathbf{J}:T^\ast
Q\to \mathfrak{g}^\ast$, the regular orbit reduced
space $(T^\ast Q)_{\mathcal{O}_\mu}=
\mathbf{J}^{-1}(\mathcal{O}_\mu)/G$ is a regular quotient symplectic
manifold with the symplectic form $\omega_{\mathcal{O}_\mu}$
uniquely characterized by the relation
\begin{equation}i_{\mathcal{O}_\mu}^\ast \cdot \omega_0=\pi_{\mathcal{O}_{\mu}}^\ast
\cdot \omega_{\mathcal{O}
_\mu}+\mathbf{J}_{\mathcal{O}_\mu}^\ast \cdot \omega_{\mathcal{O}_\mu}^+,
\label{2.7}
\end{equation}
where $\mathbf{J}_{\mathcal{O}_\mu}$ is
the restriction of the momentum map $\mathbf{J}$ to
$\mathbf{J}^{-1}(\mathcal{O}_\mu)$, that is,
$\mathbf{J}_{\mathcal{O}_\mu}=\mathbf{J}\cdot
i_{\mathcal{O}_\mu}$,
and $\omega_{\mathcal{O}_\mu}^+$ is the $+$-symplectic structure on
the orbit $\mathcal{O}_\mu$ given by
\begin{equation}\omega_{\mathcal{O}_\mu}^
+(\nu)(\xi_{\mathfrak{g}^\ast}(\nu),\eta_{\mathfrak{g}^\ast}(\nu))
=<\nu,[\xi,\eta]>,\;\; \forall\;\nu\in\mathcal{O}_\mu, \;
\xi,\eta\in \mathfrak{g}, \; \xi_{\mathfrak{g}^\ast}, \eta_{\mathfrak{g}^\ast}
\in \mathfrak{g}^*. \label{2.8}
\end{equation} The maps
$i_{\mathcal{O}_\mu}:\mathbf{J}^{-1}(\mathcal{O}_\mu)\rightarrow
T^\ast Q$ and
$\pi_{\mathcal{O}_\mu}:\mathbf{J}^{-1}(\mathcal{O}_\mu)\rightarrow
(T^\ast Q)_{\mathcal{O}_\mu}$ are natural injection and the
projection, respectively. In general case,
we maybe thought that the structure of the
symplectic orbit reduced space $((T^\ast
Q)_{\mathcal{O}_\mu},\omega_{\mathcal{O}_\mu})$ is more complex than
that of the symplectic point reduced space $((T^\ast
Q)_\mu,\omega_\mu)$, but, from Ortega and Ratiu \cite{orra04}
and the regular reduction diagram, we
know that the regular orbit reduced space $((T^\ast
Q)_{\mathcal{O}_\mu},\omega_{\mathcal{O}_\mu})$ is symplectically
diffeomorphic to the regular point reduced space $((T^*Q)_\mu,
\omega_\mu)$, and hence is also symplectically diffeomorphic to a
symplectic fiber bundle. Comparing the regular orbit reduced spaces
$((TQ)_{\mathcal{O}_\mu},\omega^L_{\mathcal{O}_\mu})$ and
$((T^\ast Q)_{\mathcal{O}_\mu},\omega_{\mathcal{O}_\mu})$ at the orbit $\mathcal{O}_\mu$,
we have that the following theorem holds.
\begin{theo}
Assume that the Lagrangian $L: TQ\rightarrow \mathbb{R}$ is
hyperregular, and the Legendre transformation $\mathcal{F}L: TQ
\rightarrow T^{\ast}Q $ is $(\Phi^T,\; \Phi^{T*})$-equivariant,
then the regular orbit reduced space
$((TQ)_{\mathcal{O}_\mu},\omega^L_{\mathcal{O}_\mu})$
of $(TQ, \omega^L)$ at the orbit $\mathcal{O}_\mu$ is
symplectically diffeomorphic to the regular orbit reduced space $((T^\ast
Q)_{\mathcal{O}_\mu},\omega_{\mathcal{O}_\mu})$
of $(T^*Q, \omega_0)$ at the orbit $\mathcal{O}_\mu$
and to the regular point reduced space $((T^*Q)_\mu, \omega_\mu)$ of $(T^*Q, \omega_0)$ at $\mu$,
and hence is also symplectically diffeomorphic to a
symplectic fiber bundle.
\end{theo}

\noindent {\bf Proof:} Since the Lagrangian $L: TQ\rightarrow \mathbb{R}$ is
hyperregular, then the Legendre transformation $\mathcal{F}L: TQ
\rightarrow T^{\ast}Q $ is a diffeomorphism. Because $\mathcal{F}L$
is $(\Phi^T,\; \Phi^{T*})$-equivariant, that is,
$\Phi^{T*}\cdot \mathcal{F}L= \mathcal{F}L \cdot \Phi^{T}$,
then we can define a map
$(\mathcal{F}L)_{\mathcal{O}_\mu}: (TQ)_{\mathcal{O}_\mu}\rightarrow (T^*Q)_{\mathcal{O}_\mu} $
given by
$(\mathcal{F}L)_{\mathcal{O}_\mu} \cdot \tau_{\mathcal{O}_\mu}
= \pi_{\mathcal{O}_\mu} \cdot \mathcal{F}L, $
and $i_{\mathcal{O}_\mu} \cdot \mathcal{F}L= \mathcal{F}L \cdot j_{\mathcal{O}_\mu} ,$
see the following commutative Diagram-3,
which is well-defined and a diffeomorphism.
$$
\begin{CD}
\mathbf{J}_L^{-1}(\mathcal{O}_\mu)\subset TQ @> \mathcal{F}L  >> \mathbf{J}^{-1}(\mathcal{O}_\mu) \subset T^*Q \\
@V \tau_{\mathcal{O}_\mu} VV @VV \pi_{\mathcal{O}_\mu}  V \\
(TQ)_{\mathcal{O}_\mu} @> (\mathcal{F}L)_{\mathcal{O}_\mu} >> (T^*Q)_{\mathcal{O}_\mu}
\end{CD}
$$
$$\mbox{Diagram-3}$$

We shall prove that
$(\mathcal{F}L)_{\mathcal{O}_\mu}$ is symplectic, that is,
$(\mathcal{F}L)_{\mathcal{O}_\mu}^* \cdot \omega_{\mathcal{O}_\mu}
= \omega^L_{\mathcal{O}_\mu}.$
In fact, from (2.7), (2.5) and (2.4), we have that
\begin{align*}
\tau_{\mathcal{O}_\mu}^* \cdot (\mathcal{F}L)_{\mathcal{O}_\mu}^* \cdot \omega_{\mathcal{O}_\mu}
& = ( (\mathcal{F}L)_{\mathcal{O}_\mu} \cdot \tau_{\mathcal{O}_\mu})^* \cdot \omega_{\mathcal{O}_\mu}
=(\pi_{\mathcal{O}_\mu} \cdot \mathcal{F}L)^* \cdot \omega_{\mathcal{O}_\mu}\\
& = (\mathcal{F}L)^* \cdot \pi_{\mathcal{O}_\mu}^* \cdot \omega_{\mathcal{O}_\mu}
= (\mathcal{F}L)^* \cdot (i_{\mathcal{O}_\mu}^* \cdot \omega_0
- \mathbf{J}_{\mathcal{O}_\mu}^\ast \cdot \omega_{\mathcal{O}_\mu}^+)\\
& =  (\mathcal{F}L)^* \cdot i_{\mathcal{O}_\mu}^\ast \cdot \omega_0
- (\mathcal{F}L)^* \cdot (\mathbf{J}_{\mathcal{O}_\mu}^\ast \cdot \omega_{\mathcal{O}_\mu}^+)\\
& = (i_{\mathcal{O}_\mu} \cdot \mathcal{F}L)^* \cdot \omega_0
- (\mathbf{J}_{\mathcal{O}_\mu} \cdot \mathcal{F}L)^* \cdot \omega_{\mathcal{O}_\mu}^{L+}\\
& = (\mathcal{F}L \cdot j_{\mathcal{O}_\mu})^* \cdot \omega_0
- (\mathbf{J}\cdot i_{\mathcal{O}_\mu} \cdot \mathcal{F}L)^* \cdot \omega_{\mathcal{O}_\mu}^{L+}\\
& = j_{\mathcal{O}_\mu}^*\cdot  (\mathcal{F}L)^* \cdot \omega_0
- (\mathbf{J}\cdot \mathcal{F}L \cdot j_{\mathcal{O}_\mu})^*\cdot \omega_{\mathcal{O}_\mu}^{L+}\\
& =  j_{\mathcal{O}_\mu}^* \cdot \omega^L
- (\mathbf{J}_L \cdot j_{\mathcal{O}_\mu})^*\cdot \omega_{\mathcal{O}_\mu}^{L+}\\
& =  j_{\mathcal{O}_\mu}^* \cdot \omega^L
- (\mathbf{J}_L )_{\mathcal{O}_\mu}^*\cdot \omega_{\mathcal{O}_\mu}^{L+}\\
& = \tau_{\mathcal{O}_\mu}^* \cdot \omega^L_{\mathcal{O}_\mu}.
\end{align*}
Notice that $\tau_{\mathcal{O}_\mu}$ is surjective, and hence
$(\mathcal{F}L)_{\mathcal{O}_\mu}^* \cdot \omega_{\mathcal{O}_\mu}
= \omega^L_{\mathcal{O}_\mu}.$
Thus, the regular orbit reduced space $((TQ)_{\mathcal{O}_\mu},\omega^L_{\mathcal{O}_\mu})$
of $(TQ, \omega^L)$ at the orbit $\mathcal{O}_\mu$ is
symplectically diffeomorphic to the regular orbit reduced space $((T^\ast
Q)_{\mathcal{O}_\mu},\omega_{\mathcal{O}_\mu})$
of $(T^*Q, \omega_0)$ at the orbit $\mathcal{O}_\mu$.
From Ortega and Ratiu \cite{orra04}
and the regular reduction diagram, we
know that the regular orbit reduced space $((T^\ast
Q)_{\mathcal{O}_\mu},\omega_{\mathcal{O}_\mu})$ at the orbit $\mathcal{O}_\mu$
is symplectically diffeomorphic to the regular point reduced space
$((T^*Q)_\mu, \omega_\mu)$ of $(T^*Q, \omega_0)$ at $\mu$,
and hence $((TQ)_{\mathcal{O}_\mu},\omega^L_{\mathcal{O}_\mu})$
is symplectically diffeomorphic to the regular point reduced space
$((T^*Q)_\mu, \omega_\mu)$ at $\mu$, and is also symplectically diffeomorphic to a
symplectic fiber bundle.
\hskip 1cm $\blacksquare$ \\

Thus, from the above discussion, we know that the regular point or regular orbit reduced
space for a regular Lagrangian system defined on a tangent bundle
may not be a tangent bundle. Considering
the completeness of the symmetric reduction, if we may define an RCL
system on a symplectic fiber bundle, then it is possible to describe
uniformly the RCL systems on $TQ$ and their regular reduced RCL
systems on the associated reduced spaces.

\section{Regular Controlled Lagrangian System and RCL-Equivalence}

In order to give a proper definition of CL system,
following the ideas in Marsden et al. \cite{mawazh10}, we first define a CL
system on $TQ$ by using the Lagrangian symplectic form, and such system is
called a regular controlled Lagrangian (RCL) system,
and then regard a regular Lagrangian system on $TQ$
as a spacial case of an RCL system without external force and
control. Thus, the set of the regular Lagrangian systems on $TQ$ is a subset
of the set of RCL systems on $TQ$. On the other hand,
since the regular reduced system of a regular Lagrangian system with symmetry
defined on the tangent bundle $TQ$ may not be a regular Lagrangian system
on a tangent bundle. So,
we can not define directly an RCL system on the tangent bundle $TQ$.
However, from Theorem 2.7 and Theorem 2.8, we know that
the regular point reduced space $((TQ)_\mu, \omega^L_\mu)$
of $(TQ, \omega^L)$ at $\mu$ is symplectically
diffeomorphic to a symplectic fiber bundle over $T(Q/G_\mu)$
with fiber to be the co-adjoint orbit $\mathcal{O}_\mu$,
and the regular orbit reduced space
$((TQ)_{\mathcal{O}_\mu},\omega^L_{\mathcal{O}_\mu})$
of $(TQ, \omega^L)$ at the orbit $\mathcal{O}_\mu$ is also symplectically
diffeomorphic to a symplectic fiber bundle.
In consequence, if we may define an RCL system on a symplectic fiber bundle,
then it is possible to describe uniformly the RCL system on $TQ$ and its
regular reduced RCL systems on the associated reduced spaces, and
we can study regular
reduction of the RCL systems with symmetries and momentum maps,
as an extension of the regular reduction theory of the regular
Lagrangian systems under regular controlled Lagrangian equivalence
conditions, and set up the regular reduction theory of the RCL system on a tangent
bundle, by using momentum map and the associated reduced Lagrangian symplectic
form and from the viewpoint of completeness of regular reduction.\\

In this section, we first define an RCL system on a symplectic
fiber bundle, then we obtain the RCL system on a tangent
bundle as a special case, by using the Legendre transformation and
the Lagrangian symplectic form on the tangent bundle of a configuration
manifold, and give a good expression of the dynamical vector field of the RCL system,
such that we can discuss RCL-equivalence. In consequence, we can study
the RCL systems with symmetries by combining with the symmetric
reduction of the regular Lagrangian systems with symmetries. For convenience, we assume
that all controls appearing in this paper are the admissible controls.\\

Let $(E,M,N,\pi,G)$ be a fiber bundle and $(E, \omega_E)$ be a
symplectic fiber bundle. If a function $L: E \rightarrow
\mathbb{R}$ is hyperregular
Lagrangian, and there is an action function $A: E \rightarrow
\mathbb{R}$ and an Euler-Lagrange vector field $\xi_L$ satisfy the equation
$\mathbf{i}_{\xi_L}\omega_E=\mathbf{d}E_L$, where $E_L: E\rightarrow
\mathbb{R}$ is an energy function given by $E_L:= A-L$. Then $(E, \omega_E, L)$ is a
regular Lagrangian system. Moreover, if considering the external
force and control, we can define a kind of regular controlled
Lagrangian (RCL) system on the symplectic fiber bundle $E$ as
follows.

\begin{defi}
(RCL System) An RCL system on $E$ is a 5-tuple $(E, \omega_E, L, F^L,
\mathcal{C}^L)$, where $(E, \omega_E, L)$ is a regular Lagrangian
system, and the function $L: E \rightarrow \mathbb{R}$ is called the (hyperregular)
Lagrangian, a fiber-preserving map $F^L: E\rightarrow E$ is called
the (external) force map, and a fiber submanifold $\mathcal{C}^L$ of
$E$ is called the control subset.
\end{defi}
Sometimes, $\mathcal{C}^L$ is also denoted the set of fiber-preserving
maps from $E$ to $\mathcal{C}^L$. When a feedback control law $u^L:
E\rightarrow \mathcal{C}^L$ is chosen, the 5-tuple $(E, \omega_E, L,
F^L, u^L)$ is a closed-loop dynamical system. In particular, when
$Q$ is a smooth manifold, and $TQ$ its tangent bundle, and $T^\ast
Q$ its cotangent bundle with a canonical symplectic form $\omega_0$,
assume that $L:TQ\rightarrow \mathbb{R}$ is a hyperregular
Lagrangian on $TQ$, and the Legendre transformation
$\mathcal{F}L: TQ \rightarrow T^{\ast}Q $ is a diffeomorphism,
then $(TQ, \omega^L )$ is a symplectic
vector bundle, where $\omega^L= \mathcal{F}L^*(\omega_0)$. If we take
that $E= TQ$, from above definition we can obtain an RCL system on
the tangent bundle $TQ$, that is, 5-tuple $(TQ, \omega^L, L, F^L,
\mathcal{C}^L)$.\\

In order to describe the dynamics of the RCL system
$(E,\omega_E,L,F^L,\mathcal{C}^L)$ with a control law $u^L:
E\rightarrow \mathcal{C}^L$, we need to give a good expression of
the dynamical vector field of the RCL system. We shall use the notations
of vertical lift maps of a vector along a fiber introduced in
Marsden et al.\cite{mawazh10}. In fact, for a smooth manifold $E$,
its tangent bundle $TE$ is a vector bundle, and for the fiber bundle
$\pi: E \rightarrow M$, we consider the tangent mapping $T\pi: TE
\rightarrow TM$ and its kernel $ker (T\pi)=\{\rho\in TE|
T\pi(\rho)=0\}$, which is a vector subbundle of $TE$. Denote by
$VE:= ker(T\pi)$, which is called a vertical bundle of $E$. Assume
that there is a metric on $E$, and we take a Levi-Civita connection
$\mathcal{A}$ on $TE$, and denote by $HE:= ker(\mathcal{A})$, which
is called a horizontal bundle of $E$, such that $TE= HE \oplus VE. $
For any $x\in M, \; a_x, b_x \in E_x, $ any tangent vector
$\rho(b_x)\in T_{b_x}E$ can be split into horizontal and vertical
parts, that is, $\rho(b_x)= \rho^h(b_x)\oplus \rho^v(b_x)$, where
$\rho^h(b_x)\in H_{b_x}E$ and $\rho^v(b_x)\in V_{b_x}E$. Let
$\gamma$ be a geodesic in $E_x$ connecting $a_x$ and $b_x$, and
denote by $\rho^v_\gamma(a_x)$ a tangent vector at $a_x$, which is a
parallel displacement of the vertical vector $\rho^v(b_x)$ along the
geodesic $\gamma$ from $b_x$ to $a_x$. Since the angle between two
vectors is invariant under a parallel displacement along a geodesic,
then $T\pi(\rho^v_\gamma(a_x))=0, $ and hence $\rho^v_\gamma(a_x)
\in V_{a_x}E. $ Now, for $a_x, b_x \in E_x $ and tangent vector
$\rho(b_x)\in T_{b_x}E$, we can define the vertical lift map of a
vector along a fiber given by
$$\mbox{vlift}: TE_x \times E_x \rightarrow TE_x; \;\; \mbox{vlift}(\rho(b_x),a_x) = \rho^v_\gamma(a_x). $$
It is easy to check from the basic fact in differential geometry
that this map does not depend on the choice of $\gamma$. If $F^L: E
\rightarrow E$ is a fiber-preserving map, for any $x\in M$, we have
that $F^L_x: E_x \rightarrow E_x$ and $TF^L_x: TE_x \rightarrow
TE_x$, then for any $a_x \in E_x$ and $\rho\in TE_x$, the vertical
lift of $\rho$ under the action of $F^L$ along a fiber is defined by
$$(\mbox{vlift}(F^L_x)\rho)(a_x)=\mbox{vlift}((TF^L_x\rho)(F^L_x(a_x)), a_x)= (TF^L_x\rho)^v_\gamma(a_x), $$
where $\gamma$ is a geodesic in $E_x$ connecting $F^L_x(a_x)$ and
$a_x$.\\

In particular, when $\pi: E \rightarrow M$ is a vector bundle, for
any $x\in M$, the fiber $E_x=\pi^{-1}(x)$ is a vector space. In this
case, we can choose the geodesic $\gamma$ to be a straight line, and
the vertical vector is invariant under a parallel displacement along
a straight line, that is, $\rho^v_\gamma(a_x)= \rho^v(b_x).$
Moreover, when $E= TQ$, by using the local trivialization of $TTQ$,
we have that $TTQ\cong TQ \times TQ$ (locally). Because of $\tau_Q: TQ
\rightarrow Q$, and $T\tau_Q: TTQ \rightarrow TQ$, then in this case,
for any $v_x, \; w_x \in T_x Q, \; x\in Q, $ we know that $(0, w_x)
\in V_{w_x}T_x Q, $ and hence we can get that
$$ \mbox{vlift}((0, w_x)(w_x), v_x) = (0, w_x)(v_x)
= \left.\frac{\mathrm{d}}{\mathrm{d}s}\right|_{s=0}(v_x+s w_x),
$$ which coincides with the definition of vertical lift map
along fiber in Marsden and Ratiu \cite{mara99}.\\

For a given RCL System $(TQ, \omega^L, L, F^L, \mathcal{C}^L)$, the
dynamical vector field of the associated regular Lagrangian system
$(TQ, \omega^L, L) $ is the Euler-Lagrange vector field $\xi_L$, such that
$\mathbf{i}_{\xi_L}\omega^L=\mathbf{d}E_L$. If considering
the external force $F^L: TQ \rightarrow TQ, $ by using the above
notation of vertical lift map of a vector along a fiber, the
change of $\xi_L$ under the action of $F^L$ is that
$$\mbox{vlift}(F^L)\xi_L(v_x)= \mbox{vlift}((TF^L\xi_L)(F^L(v_x)), v_x)= (TF^L\xi_L)^v_\gamma(v_x),$$
where $v_x \in T_x Q, \; x\in Q $ and the geodesic $\gamma$ is a straight line in
$T_x Q$ connecting $F^L_x(v_x)$ and $v_x$. In the same way, when a
feedback control law $u^L: TQ \rightarrow \mathcal{C}^L$ is chosen,
the change of $\xi_L$ under the action of $u^L$ is that
$$\mbox{vlift}(u^L)\xi_L(v_x)= \mbox{vlift}((Tu^L\xi_L)(u^L(v_x)), v_x)= (Tu^L\xi_L)^v_\gamma(v_x).$$
In consequence, we can give an expression of the dynamical vector
field of the RCL system as follows.
\begin{theo}
The dynamical vector field of an RCL system
$(TQ,\omega^L,L,F^L,\mathcal{C}^L)$ with a control law $u^L$ is the
synthetic of the Euler-Lagrange vector field $\xi_L$ and its changes under
the actions of the external force $F^L$ and control $u^L$, that is,
$$\xi_{(TQ,\omega^L,L,F^L,u^L)}(v_x)= \xi_L(v_x)+ \textnormal{vlift}(F^L)\xi_L(v_x)+ \textnormal{vlift}(u^L)\xi_L(v_x),$$
for any $v_x \in T_x Q, \; x\in Q $. For convenience, it is simply
written as
\begin{equation}\xi_{(TQ,\omega^L,L,F^L,u^L)}
=\xi_L+\textnormal{vlift}(F^L)+\textnormal{vlift}(u^L).\label{3.1}\end{equation}
\end{theo}
Where $\textnormal{vlift}(F^L)= \textnormal{vlift}(F^L)\xi_L,$
and $\textnormal{vlift}(u^L)= \textnormal{vlift}(u^L)\xi_L$ are
the changes of $\xi_L$ under the actions of $F^L$ and $u^L$.
We also denote that $\textnormal{vlift}(\mathcal{C}^L)=
\bigcup\{\textnormal{vlift}(u^L)\xi_L | \; u^L\in \mathcal{C}^L\}$.
It is worthy of noting that, in order to deduce and calculate
easily, we always use the simple expression of dynamical vector
field $\xi_{(TQ,\omega^L,L,F^L,u^L)}$. Moreover, we also use the
simple expressions for $R_p$-reduced vector field $\xi_{((TQ)_\mu,
\omega^L_\mu, l_\mu, f^L_\mu, u^L_\mu)}$ and $R_o$-reduced vector
field $\xi_{((TQ)_{\mathcal{O}_\mu}, \omega^L_{\mathcal{O}_\mu},
l_{\mathcal{O}_\mu},f^L_{\mathcal{O}_\mu},u^L_{\mathcal{O}_\mu})}$
in Section 4 and Section 5.\\

From the expression (3.1) of the dynamical vector
field of the RCL system, we know that under the actions of the external force $F^L$
and control $u^L$, in general, the dynamical vector
field may not be an Euler-Lagrange vector field, and hence the RCL system may not
be yet a regular Lagrangian system. However,
it is a dynamical system closed relative to a
regular Lagrangian system, and it can be explored and studied by extending
the methods for external force and control
in the study of the regular Lagrangian system.
In particular, it is worthy of noting that
the energy $E_L$ is conserved along the flow of the
Euler-Lagrange vector field $\xi_L$, if $\xi_L$ satisfies the second order
equation $T\tau_Q \circ \xi_L = id_{TQ}.$ Note that
$T\tau_Q \cdot \textnormal{vlift}(F^L) = T\tau_Q \cdot \textnormal{vlift}(u^L) =0,$
then from the expression (3.1) we have that
$T\tau_Q \circ \xi_{(TQ,\omega^L,L,F^L,u^L)} = id_{TQ}, $ that is,
the dynamical vector field of the RCL system satisfies always the second order equation.\\

On the other hand, for two given regular Lagrangian systems
$(TQ_i,\omega^L_i,L_i),$ $ i= 1,2,$ we say them to be
equivalent, if there exists a
diffeomorphism $\varphi: Q_1\rightarrow Q_2$, such that
their Euler-Lagrange vector fields $\xi_{L_i}, \; i=1,2 $ satisfy
the condition $\xi_{L_2}\cdot T\varphi
=T(T\varphi) \cdot \xi_{L_1}$, where the map
$T\varphi: TQ_1 \rightarrow TQ_2$
is the tangent map of $\varphi$, and the
map $T(T\varphi): TTQ_1 \rightarrow TTQ_2$ is the
tangent map of $T\varphi$. It is easy to see that the condition
$\xi_{L_2}\cdot T\varphi =T(T\varphi) \cdot \xi_{L_1}$ is equivalent the fact that
the map $T\varphi: TQ_1\rightarrow TQ_2$
is symplectic with respect to their Lagrangian symplectic forms
$\omega^{L}_{i}$ on $TQ_i, \; i=1,2.$ \\

For two given RCL systems $(TQ_i,\omega^{L}_i,L_i,F^{L}_i,
\mathcal{C}^{L}_i ),$ $ i= 1,2,$ we also want to define
their equivalence, that is, to look for a diffeomorphism
$\varphi: Q_1\rightarrow Q_2$, such that
$\xi_{(TQ_2,\omega^{L}_2,L_2,F^{L}_2,\mathcal{C}^{L}_2)}\cdot T\varphi
=T(T\varphi)\cdot \xi_{(TQ_1,\omega^{L}_1,L_1,F^{L}_1,\mathcal{C}^{L}_1)}$. But,
it is worthy of noting that,
when an RCL system is given, the force map $F^L: TQ \rightarrow TQ$
is determined, but the feedback control law $u^L: TQ\rightarrow
\mathcal{C}^L$ could be chosen. In order to emphasize explicitly
the impact of external force
and control in the study of the RCL systems, by using
the above expression (3.1) of the dynamical vector field of the RCL system,
we can describe the feedback control law how
to modify the structure of the RCL system, and the regular controlled
Lagrangian matching conditions and RCL-equivalence are induced as
follows.
\begin{defi}
(RCL-equivalence) Suppose that we have two RCL systems
$(TQ_i,\omega^{L}_i,L_i,F^{L}_i,\mathcal{C}^{L}_i) $, $i= 1,2,$ we say
them to be RCL-equivalent, or simply,
$(TQ_1,\omega^{L}_1,L_1,F^{L}_1,\mathcal{C}^{L}_1)\stackrel{RCL}{\sim}\\
(TQ_2,\omega^{L}_2,L_2,F^{L}_2,\mathcal{C}^{L}_2)$,
if there exists a diffeomorphism $\varphi: Q_1\rightarrow Q_2$, such
that the following regular controlled Lagrangian matching conditions hold:

\noindent {\bf RCL-1:} The control subsets $\mathcal{C}^{L}_i, \; i=1,2$ satisfy
the condition $\mathcal{C}^{L}_2=T\varphi (\mathcal{C}^{L}_1),$ where the map
$T\varphi: TQ_1 \rightarrow TQ_2$ is tangent map of $\varphi$.

\noindent {\bf RCL-2:} For each control law $u^{L}_1:
TQ_1 \rightarrow \mathcal{C}^{L}_1, $ there exists the control law $u^{L}_2:
TQ_2 \rightarrow \mathcal{C}^{L}_2, $  such that the two
closed-loop dynamical systems produce the same dynamical vector fields, that is,
$\xi_{(TQ_2,\omega^{L}_2,L_2,F^{L}_2,u^{L}_2)}\cdot T\varphi
=T(T\varphi)\cdot \xi_{(TQ_1,\omega^{L}_1,L_1,F^{L}_1,u^{L}_1)}$,
where the map $T(T\varphi): TTQ_1 \rightarrow TTQ_2$
is the tangent map of $T\varphi$.
\end{defi}

From the expression (3.1) of the dynamical vector field of the RCL system
and the condition $\xi_{(TQ_2,\omega^{L}_2,L_2,F^{L}_2,u^{L}_2)}\cdot T\varphi
=T(T\varphi)\cdot \xi_{(TQ_1,\omega^{L}_1,L_1,F^{L}_1,u^{L}_1)}$, we have that
$$
(\xi_{L_2}+\textnormal{vlift}(F^{L}_2)\xi_{L_2}+\textnormal{vlift}(u^{L}_2)\xi_{L_2} )\cdot T\varphi
= T(T\varphi)\cdot [\xi_{L_1}+\textnormal{vlift}(F^{L}_1)\xi_{L_1}+\textnormal{vlift}(u^{L}_1)\xi_{L_1}].
$$
By using the notation of vertical lift map of a vector along a fiber,
for $v_x \in T_x Q_1, \; x \in Q_1$, we have that
\begin{align*}
&T(T\varphi)\cdot \textnormal{vlift}(F^{L}_1)\xi_{L_1}(v_x) \\
&=T(T\varphi)\cdot \textnormal{vlift}((TF^{L}_1\cdot \xi_{L_1}) (F^{L}_1(v_x)), v_x)\\
&=\textnormal{vlift}((T(T\varphi)\cdot TF^{L}_1\cdot T(T\varphi^{-1})\cdot \xi_{L_1})
(T\varphi \cdot F^{L}_1 \cdot T\varphi^{-1}\cdot(T\varphi \cdot v_x)), T\varphi \cdot v_x)\\
&=\textnormal{vlift}((T((T\varphi) \cdot F^{L}_1 \cdot T\varphi^{-1}) \cdot \xi_{L_1})
(T\varphi \cdot F^{L}_1 \cdot T\varphi^{-1}(T\varphi \cdot v_x)), T\varphi \cdot v_x)\\
&=\textnormal{vlift}(T\varphi \cdot
F^{L}_1 \cdot T\varphi^{-1})\cdot \xi_{L_1}(T\varphi \cdot v_x),
\end{align*}
where the map $T\varphi^{-1}: TQ_2 \rightarrow TQ_1$.
In the same way, we have that
$T(T\varphi)\cdot \textnormal{vlift}(u^{L}_1)\xi_{L_1}=\textnormal{vlift}(T\varphi
\cdot u^{L}_1 \cdot T\varphi^{-1})\cdot \xi_{L_1}\cdot T\varphi $.
Thus, the explicit relation between the two control laws
$u^{L}_i: TQ_i \rightarrow \mathcal{C}^L_i, \; i=1,2$ in {\bf RCL-2} is given by
\begin{align}
& (\textnormal{vlift}(u^{L}_2)
-\textnormal{vlift}(T\varphi \cdot u^{L}_1 \cdot T\varphi^{-1}))\cdot T\varphi \nonumber \\
& = -\xi_{L_2} \cdot T\varphi +T(T\varphi) (\xi_{L_1})+
(-\textnormal{vlift}(F^{L}_2)+\textnormal{vlift}(T\varphi \cdot F^{L}_1
\cdot T\varphi^{-1}))\cdot T\varphi.
\label{3.2}\end{align}

From the above relation we know that, when two RCL systems
$(TQ_i,\omega^{L}_i,L_i, F^{L}_i,\mathcal{C}^{L}_i) $, $i= 1,2,$ are
RCL-equivalent with respect to $T\varphi$, the corresponding regular Lagrangian systems
$(TQ_i,\omega^{L}_i,L_i),$ $ i= 1,2,$ may not be equivalent with respect to $T\varphi$. If
two corresponding regular Lagrangian systems are also equivalent with respect to $T\varphi$,
then the control laws $u^{L}_i: TQ_i \rightarrow \mathcal{C}^L_i, \; i=1,2$ and the
external forces $F^{L}_i: TQ_i \rightarrow TQ_i, \; i=1,2$ in {\bf RCL-2}
must satisfy the following condition
\begin{equation}
\textnormal{vlift}(u^{L}_2) -\textnormal{vlift}(T\varphi \cdot u^{L}_1 \cdot T\varphi^{-1})
= -\textnormal{vlift}(F^{L}_2)+\textnormal{vlift}(T\varphi \cdot F^{L}_1 \cdot T\varphi^{-1}).
\label{3.3}
\end{equation}

In the following we shall introduce the regular point and regular
orbit reducible RCL systems with symmetries, and show a variety of
relationships of their regular reducible RCL-equivalences.

\section{Regular Point Reduction of the RCL System}

We know that, when the external force and control of an RCL
system $(TQ,\omega^L,L,F^L,\mathcal{C}^L)$ are both zeros,
that is, $F^L=0,$ and $\mathcal{C}^L= \emptyset,$ in this case the RCL system
is just a regular Lagrangian system $(TQ,\omega^L,L)$.
Thus, we can regard a regular Lagrangian system on $TQ$
as a spacial case of the RCL system without external force and
control. In consequence, the set of regular Lagrangian systems with symmetries
on $TQ$ is a subset of the set of RCL systems with symmetries on $TQ$.
If we first admit the regular point reduction of a regular Lagrangian system
with symmetry, then we may study the regular point
reduction of an RCL system with symmetry,
as an extension of the regular point reduction of
a regular Lagrangian system under the regular controlled Lagrangian equivalence
conditions. In order to do these, in this section
we consider the RCL system with symmetry and momentum map,
and first give the regular point reducible RCL system
and the RpCL-equivalence, then prove the regular point reduction theorems
for the RCL system and regular Lagrangian system.\\

We know that, if an RCL system with symmetry and momentum map is
regular point reducible, then the associated regular
Lagrangian system must be regular point reducible.
Thus, from Definition 2.4 and Theorem 2.7,
if the Legendre transformation $\mathcal{F}L: TQ
\rightarrow T^{\ast}Q $ is $(\Phi^T,\; \Phi^{T*})$-equivariant,
then we can introduce
a kind of regular point reducible RCL systems as follows.
\begin{defi}
(Regular Point Reducible RCL System) A 6-tuple $(TQ, G, \omega^L, L,
F^L, \mathcal{C}^L)$, where the hyperregular Lagrangian $L: TQ \rightarrow
\mathbb{R}$, the fiber-preserving map $F^L: TQ \rightarrow TQ$ and
the fiber submanifold $\mathcal{C}^L$ of $TQ$ are all $G$-invariant,
is called a {\bf regular point reducible RCL system}, if
the Legendre transformation $\mathcal{F}L: TQ
\rightarrow T^{\ast}Q $ is $(\Phi^T,\; \Phi^{T*})$-equivariant,
and there exists a
point $\mu\in\mathfrak{g}^\ast$, which is a regular value of the
momentum map $\mathbf{J}_L$, such that the regular point reduced
system, that is, the 5-tuple $((TQ)_\mu,
\omega^L_\mu,l_\mu,f^L_\mu,\mathcal{C}^L_\mu)$, where $(T
Q)_\mu=\mathbf{J}^{-1}_L(\mu)/G_\mu$, $\tau_\mu^\ast
\omega^L_\mu=j_\mu^\ast\omega^L$, $l_\mu\cdot \tau_\mu= L\cdot
j_\mu$, $F^L(\mathbf{J}^{-1}_L(\mu))\subset \mathbf{J}^{-1}_L(\mu)
$, $f^L_\mu\cdot \tau_\mu=\tau_\mu \cdot F^L\cdot j_\mu$,
$\mathcal{C}^L \cap \mathbf{J}^{-1}_L(\mu)\neq \emptyset $,
$\mathcal{C}^L_\mu=\tau_\mu(\mathcal{C}^L\cap
\mathbf{J}^{-1}_L(\mu))$, is an RCL system, which is simply written
as $R_p$-reduced RCL system. Where $((TQ)_\mu,\omega^L_\mu)$ is the
$R_p$-reduced space, the function $l_\mu: (TQ)_\mu \rightarrow
\mathbb{R}$ is called the $R_p$-reduced Lagrangian, the fiber-preserving
map $f^L_\mu:(TQ)_\mu\rightarrow (TQ)_\mu$ is called the $R_p$-reduced
(external) force map, $\mathcal{C}^L_\mu$ is a fiber submanifold of
\;$(TQ)_\mu$ and is called the $R_p$-reduced control subset.
\end{defi}

It is worthy of noting that for the regular point reducible RCL system
$(TQ,G,\omega^L,L,F^L,\mathcal{C}^L)$, the $G$-invariant external
force map $F^L: TQ \rightarrow TQ $ has to satisfy the conditions
$F^L(\mathbf{J}^{-1}_L(\mu))\subset \mathbf{J}^{-1}_L(\mu), $ and
$f^L_\mu\cdot \tau_\mu=\tau_\mu \cdot F^L\cdot j_\mu, $ such that we
can define the $R_p$-reduced external force map $f^L_\mu:(TQ)_\mu
\rightarrow (TQ)_\mu. $ The condition $\mathcal{C}^L \cap
\mathbf{J}^{-1}_L(\mu)\neq \emptyset $ in above definition makes
that the $G$-invariant control subset $\mathcal{C}^L\cap
\mathbf{J}^{-1}_L(\mu)$ can be reduced and the $R_p$-reduced control
subset is $\mathcal{C}^L_\mu= \tau_\mu(\mathcal{C}^L\cap
\mathbf{J}^{-1}_L(\mu))$.\\

Assume that the dynamical vector field $\xi_{(TQ,G,\omega^L,L,F^L,u^L)}$ of
a given regular point reducible RCL system $(TQ,G,\omega^L, L, F^L,
\mathcal{C}^L)$ with a control law $u^L \in \mathcal{C}^L$ can be expressed by
\begin{equation}\xi_{(TQ,G,\omega^L,L,F^L,u^L)}
=\xi_L+\textnormal{vlift}(F^L)+\textnormal{vlift}(u^L).
\label{4.1}\end{equation}
Then, for the regular point reducible RCL system we can also
introduce the regular point reducible controlled Lagrangian
equivalence (RpCL-equivalence) as follows.
\begin{defi}
(RpCL-equivalence)
Suppose that we have two regular point reducible RCL systems $(TQ_i,
G_i,\omega^L_i,L_i, F^L_i, \mathcal{C}^L_i)$, $i=1,2$, we say them
to be RpCL-equivalent, or simply,\\ $(TQ_1,
G_1,\omega^L_1,L_1,F^L_1,\mathcal{C}^L_1)\stackrel{RpCL}{\sim}(T
Q_2,G_2,\omega^L_2,L_2,F^L_2,\mathcal{C}^L_2)$, if there exists a
diffeomorphism $\varphi:Q_1\rightarrow Q_2$ such that the following
regular point reducible controlled Lagrangian matching conditions hold:\\

\noindent {\bf RpCL-1:} For $\mu_i\in \mathfrak{g}^\ast_i $, the
regular reducible points of the RCL systems $(TQ_i, G_i,\omega^L_i, L_i,
F^L_i, \mathcal{C}^L_i) $, $i=1,2$, the map
$(T\varphi)_\mu=j_{\mu_2}^{-1}\cdot T\varphi \cdot j_{\mu_1}:
(\mathbf{J}_L)_1^{-1}(\mu_1)\rightarrow
(\mathbf{J}_L)_2^{-1}(\mu_2)$ is
$(G_{1\mu_1},G_{2\mu_2})$-equivariant and $\mathcal{C}^L_2\cap
(\mathbf{J}_L)_2^{-1}(\mu_2)=(T\varphi)_\mu (\mathcal{C}^L_1\cap
(\mathbf{J}_L)_1^{-1}(\mu_1))$, where $\mu=(\mu_1, \mu_2)$, and
denote by $j_{\mu_2}^{-1}(S)$ the pre-image of a subset $S\subset
TQ_2$ for the map $j_{\mu_2}:
(\mathbf{J}_L)_2^{-1}(\mu_2)\rightarrow TQ_2$.

\noindent {\bf RpCL-2:}
For each control law $u^{L}_1:
TQ_1 \rightarrow \mathcal{C}^{L}_1, $ there exists the control law $u^{L}_2:
TQ_2 \rightarrow \mathcal{C}^{L}_2, $  such that the two
closed-loop dynamical systems produce the same dynamical vector fields, that is,
$\xi_{(TQ_2,G_2,\omega^{L}_2,L_2,F^{L}_2,u^{L}_2)}\cdot T\varphi
=T(T\varphi)\cdot \xi_{(TQ_1,G_1,\omega^{L}_1,L_1,F^{L}_1,u^{L}_1)}$.
\end{defi}

It is worthy of noting that for the regular point reducible RCL
system, the induced equivalent map $T\varphi$ also keeps the equivariance
of $G$-action at the regular point. If a feedback control law
$u^L: TQ \rightarrow \mathcal{C}^L $ is chosen, and
$u^L \in \mathcal{C}^L \cap \mathbf{J}^{-1}_L(\mu),$
and $\mathcal{C}^L \cap \mathbf{J}^{-1}_L(\mu)\neq \emptyset $,
then the $R_p$-reduced control law $u^L_\mu:(TQ)_\mu\rightarrow
\mathcal{C}^L_\mu =\tau_\mu(\mathcal{C}^L\cap
\mathbf{J}^{-1}_L(\mu))$, and $u^L_{\mu}\cdot\tau_{\mu}=\tau_{\mu}\cdot
u^L\cdot j_{\mu}$.
The $R_p$-reduced RCL system $((TQ)_\mu,
\omega^L_\mu, l_\mu, f^L_\mu, u^L_\mu)$ is a closed-loop regular
dynamical system with the $R_p$-reduced control law $u^L_\mu$. Assume that its dynamical vector
field $\xi_{((TQ)_\mu, \omega^L_\mu, l_\mu, f^L_\mu, u^L_\mu)}$ can
be expressed by
\begin{equation}\xi_{((TQ)_\mu, \omega^L_\mu, l_\mu, f^L_\mu, u^L_\mu)}
=\xi_{l_\mu}+\textnormal{vlift}(f^L_\mu)+\textnormal{vlift}(u^L_\mu),
\label{4.2}\end{equation}
where $\xi_{l_\mu}$ is the $R_p$-reduced
Euler-Lagrange vector field, $\textnormal{vlift}(f^L_\mu)=
\textnormal{vlift}(f^L_\mu)\xi_{l_\mu}$,
$\textnormal{vlift}(u^L_\mu)=
\textnormal{vlift}(u^L_\mu)\xi_{l_\mu}$
are the changes of $\xi_{l_\mu}$ under the actions of the
$R_p$-reduced external force $f^L_\mu$
and the $R_p$-reduced control law $u^L_\mu$,
and the dynamical vector fields of the RCL system and the $R_p$-reduced
RCL system satisfy the condition
\begin{equation}\xi_{((TQ)_\mu, \omega^L_\mu, l_\mu, f^L_\mu,
u^L_\mu)}\cdot \tau_\mu= T\tau_\mu\cdot
\xi_{(TQ,G,\omega^L,L,F^L,u^L)}\cdot j_\mu,
\label{4.3}
\end{equation}
see Marsden et al. \cite{mawazh10}, Wang \cite{wa18}.
Then we can obtain the following regular point reduction theorem for
the RCL system, which explains the relationship between the
RpCL-equivalence for the regular point reducible RCL system with
symmetry and the RCL-equivalence for the associated $R_p$-reduced RCL
system.
\begin{theo}
Two regular point reducible RCL systems $(TQ_i, G_i, \omega^L_i,
L_i, F^L_i, \mathcal{C}^L_i)$, $i=1,2,$ are RpCL-equivalent if and
only if the associated $R_p$-reduced RCL systems
$((TQ_i)_{\mu_i},\omega^L_{i\mu_i},l_{i\mu_i},f^L_{i\mu_i},\\
\mathcal{C}^L_{i\mu_i}), i=1,2,$ are RCL-equivalent.
\end{theo}

\noindent {\bf Proof:} If $(TQ_1, G_1, \omega^L_1, L_1, F^L_1,
\mathcal{C}^L_1)\stackrel{RpCL}{\sim}(TQ_2, G_2, \omega^L_2, L_2, F^L_2,
\mathcal{C}^L_2)$, then there exists a diffeomorphism $\varphi:Q_1\rightarrow
Q_2$ such that for $\mu_i\in \mathfrak{g}^\ast_i, i=1,2$,
$(T\varphi)_\mu=j_{\mu_2}^{-1}\cdot T\varphi \cdot
j_{\mu_1}:(\mathbf{J}_L)_1^{-1}(\mu_1)\rightarrow
(\mathbf{J}_L)_2^{-1}(\mu_2)$ is $(G_{1\mu_1},G_{2\mu_2})$-equivariant,
and $\mathcal{C}^L_2\cap (\mathbf{J}_L)_2^{-1}(\mu_2)=(T\varphi)_\mu (\mathcal{C}^L_1\cap
(\mathbf{J}_L)_1^{-1}(\mu_1))$ and RpCL-2 holds. From the following
commutative Diagram-4:
\[
\begin{CD}
TQ_1 @<j_{\mu_1}<<( \mathbf{J}_L)_1^{-1}(\mu_1) @>\tau_{\mu_1}>> (TQ_1)_{\mu_1}\\
@V T\varphi VV @V (T\varphi)_\mu VV @V (T\varphi)_{\mu/G}VV\\
TQ_2 @<j_{\mu_2}<< (\mathbf{J}_L)_2^{-1}(\mu_2)
@>\tau_{\mu_2}>>(TQ_2)_{\mu_2}
\end{CD}
\]
$$\mbox{Diagram-4}$$
we can define a map $(T\varphi)_{\mu/G}:(T
Q_1)_{\mu_1}\rightarrow (TQ_2)_{\mu_2}$ such that
$(T\varphi)_{\mu/G}\cdot
\tau_{\mu_1}=\tau_{\mu_2}\cdot (T\varphi)_\mu$. Because
$(T\varphi)_\mu: (\mathbf{J}_L)_1^{-1}(\mu_1)\rightarrow
(\mathbf{J}_L)_2^{-1}(\mu_2)$ is $(G_{1\mu_1},G_{2\mu_2})$-equivariant,
$(T\varphi)_{\mu/G}$ is well-defined. We shall show that
$\mathcal{C}^L_{2\mu_2}=(T\varphi)_{\mu/G}(\mathcal{C}^L_{1\mu_1})$. In fact,
because $(TQ_i, G_i, \omega^L_i,
L_i, F^L_i,\mathcal{C}^L_i)$, $i=1,2,$ are regular point reducible RCL systems,
then $\mathcal{C}^L_i\cap (\mathbf{J}_L)_i^{-1}(\mu_i)\neq
\emptyset $ and $\mathcal{C}^L_{i\mu_i}=\tau_{\mu_i}(\mathcal{C}^L_i\cap
(\mathbf{J}_L)_i^{-1}(\mu_i)),\; i=1,2$. From $\mathcal{C}^L_2\cap
(\mathbf{J}_L)_2^{-1}(\mu_2)=(T\varphi)_\mu (\mathcal{C}^L_1\cap
(\mathbf{J}_L)_1^{-1}(\mu_1))$, we have that
\begin{align*}
\mathcal{C}^L_{2\mu_2} & =\tau_{\mu_2}(\mathcal{C}^L_2\cap (\mathbf{J}_L)_2^{-1}(\mu_2))
=\tau_{\mu_2}\cdot (T\varphi)_\mu (\mathcal{C}^L_1\cap
(\mathbf{J}_L)_1^{-1}(\mu_1))\\ & =(T\varphi)_{\mu/G}\cdot
\tau_{\mu_1}(\mathcal{C}^L_1\cap (\mathbf{J}_L)_1^{-1}(\mu_1))
=(T\varphi)_{\mu/G}(\mathcal{C}^L_{1\mu_1}).
\end{align*}
Thus, the condition RCL-1 holds. On the other hands,
for the $R_p$-reduced control law $u^L_{1\mu_1}:(T
Q_1)_{\mu_1}\rightarrow \mathcal{C}^L_{1\mu_1}$, we have the control law
$u^L_1: TQ_1 \rightarrow \mathcal{C}^L_1, $, such that
$u^L_{1\mu_1}\cdot\tau_{\mu_1}=\tau_{\mu_1}\cdot u^L_1\cdot
j_{\mu_1}$. From the condition RpCL-2 we know that
there exists the control law $u^L_2:
TQ_2 \rightarrow \mathcal{C}^L_2, $  such that
$\xi_{(TQ_2,G_2,\omega^L_2,L_2,F^L_2,u^L_2)}\cdot T\varphi
=T(T\varphi)\cdot \xi_{(TQ_1,G_1,\omega^L_1,L_1,F^L_1,u^L_1)}$.
But, for the control law $u^L_2: TQ_2 \rightarrow \mathcal{C}^L_2, $ we have
the $R_p$-reduced control law $u^L_{2\mu_2}:(T
Q_2)_{\mu_2}\rightarrow \mathcal{C}^L_{2\mu_2}$, such that
$u^L_{2\mu_2}\cdot\tau_{\mu_2}=\tau_{\mu_2}\cdot u^L_2\cdot
j_{\mu_2}$. Note that for $i=1,2$,
from (4.3), we have that
\begin{align}
\xi_{((TQ_i)_{\mu_i}, \omega^L_{i\mu_i}, l_{i\mu_i},
f^L_{i\mu_i},u^L_{i\mu_i})}\cdot \tau_{\mu_i}=T\tau_{\mu_i}\cdot
\xi_{(TQ_i, G_i, \omega^L_i, L_i, F^L_i, u^L_i)}\cdot j_{\mu_i},
\end{align}
and from the commutative Diagram-4,
$(T\varphi)_{\mu/G}\cdot
\tau_{\mu_1}=\tau_{\mu_2}\cdot (T\varphi)_\mu$ and
$j_{\mu_2}\cdot (T\varphi)_\mu = (T\varphi) \cdot j_{\mu_1}$,
then we have that
\begin{align*}
& \xi_{((TQ_2)_{\mu_2}, \omega^L_{2\mu_2}, l_{2\mu_2},
f^L_{2\mu_2},u^L_{2\mu_2})}\cdot (T\varphi)_{\mu/G}\cdot \tau_{\mu_1}
= \xi_{((TQ_2)_{\mu_2}, \omega^L_{2\mu_2}, l_{2\mu_2},
f^L_{2\mu_2},u^L_{2\mu_2})}\cdot \tau_{\mu_2}\cdot (T\varphi)_\mu \\
& =T\tau_{\mu_2}\cdot \xi_{(TQ_2, G_2,
\omega^L_2, L_2, F^L_2, u^L_2)}\cdot j_{\mu_2}\cdot (T\varphi)_\mu
=T\tau_{\mu_2}\cdot \xi_{(TQ_2, G_2,
\omega^L_2, L_2, F^L_2, u^L_2)}\cdot (T\varphi) \cdot j_{\mu_1}\\
& = T\tau_{\mu_2}\cdot T(T\varphi) \cdot
\xi_{(TQ_1,G_1,\omega^L_1,L_1,F^L_1,u^L_1)}\cdot j_{\mu_1}
= T(\tau_{\mu_2}\cdot (T\varphi)_\mu)\cdot
\xi_{(TQ_1,G_1,\omega^L_1,L_1,F^L_1,u^L_1)}\cdot j_{\mu_1}\\
& = T((T\varphi)_{\mu/G}\cdot \tau_{\mu_1})\cdot
\xi_{(TQ_1,G_1,\omega^L_1,L_1,F^L_1,u^L_1)}\cdot j_{\mu_1}
=T((T\varphi)_{\mu/G})\cdot T\tau_{\mu_1}\cdot
\xi_{(TQ_1,G_1,\omega^L_1,L_1,F^L_1,u^L_1)}\cdot j_{\mu_1}\\
& = T((T\varphi)_{\mu/G})\cdot
\xi_{((TQ_1)_{\mu_1}, \omega^L_{1\mu_1}, l_{1\mu_1},
f^L_{1\mu_1},u^L_{1\mu_1})}\cdot \tau_{\mu_1}.
\end{align*}
Since $\tau_{\mu_1}: (\mathbf{J}_L)_1^{-1}(\mu_1) \rightarrow
(TQ_1)_{\mu_1}$ is surjective, thus,
\begin{align}
\xi_{((TQ_2)_{\mu_2}, \omega^L_{2\mu_2}, l_{2\mu_2},
f^L_{2\mu_2},u^L_{2\mu_2})}\cdot (T\varphi)_{\mu/G}
=T((T\varphi)_{\mu/G})\cdot
\xi_{((TQ_1)_{\mu_1}, \omega^L_{1\mu_1}, l_{1\mu_1},
f^L_{1\mu_1},u^L_{1\mu_1})},
\end{align}
that is, the condition RCL-2 holds.
So, the $R_p$-reduced RCL systems $((T
Q_i)_{\mu_i},\omega^L_{i\mu_i},l_{i\mu_i}, f^L_{i\mu_i},
\mathcal{C}^L_{i\mu_i}),$ $i=1,2,$ are RCL-equivalent.\\

Conversely, assume that the $R_p$-reduced RCL systems $((T
Q_i)_{\mu_i},\omega^L_{i\mu_i},l_{i\mu_i},f^L_{i\mu_i},\mathcal{C}^L_{i\mu_i})$,
$i=1,2,$ are RCL-equivalent, then there exists a diffeomorphism
$(T\varphi)_{\mu/G}: (TQ_1)_{\mu_1}\rightarrow (T
Q_2)_{\mu_2}$, such that
$\mathcal{C}^L_{2\mu_2}=(T\varphi)_{\mu/G}(\mathcal{C}^L_{1\mu_1}),\;
\mu_i\in\mathfrak{g}_i^\ast, \; i=1,2$ and
for the $R_p$-reduced control law $u^L_{1\mu_1}:(T
Q_1)_{\mu_1}\rightarrow \mathcal{C}^L_{1\mu_1}$, there exists
the $R_p$-reduced control law $u^L_{2\mu_2}:(T
Q_2)_{\mu_2}\rightarrow \mathcal{C}^L_{2\mu_2}$, such that (4.5) holds.
Then from commutative Diagram-4, we can
define a map $(T\varphi)_\mu: (\mathbf{J}_L)^{-1}_1(\mu_1)\rightarrow
(\mathbf{J}_L)^{-1}_2(\mu_2)$ such that $\tau_{\mu_2}\cdot
(T\varphi)_\mu=(T\varphi)_{\mu/G}\cdot \tau_{\mu_1}, $ and the
map $T\varphi: TQ_1\rightarrow TQ_2 $ such that
$T\varphi \cdot j_{\mu_1}=j_{\mu_2}\cdot (T\varphi)_{\mu}, $
as well as a diffeomorphism $\varphi:
Q_1\rightarrow Q_2, $ whose tangent lift is just
$ T\varphi: TQ_1 \rightarrow TQ_2$. Moreover, for above definition of
$(T\varphi)_\mu$, we know that $(T\varphi)_\mu$ is
$(G_{1\mu_1},G_{2\mu_2})$-equivariant. In fact, for any $z_i\in
(\mathbf{J}_L)_i^{-1}(\mu_i)$, $g_i\in G_{i\mu_i}$, $i=1,2$ such that
$z_2=(T\varphi)_\mu (z_1)$, $[z_2]=(T\varphi)_{\mu/G}[z_1]$, then
we have that
\begin{align*}
(T\varphi)_\mu (\Phi_{1g_1}(z_1))& = \tau_{\mu_2}^{-1}\cdot
\tau_{\mu_2}\cdot (T\varphi)_\mu (\Phi_{1g_1}(z_1))
= \tau_{\mu_2}^{-1}\cdot \tau_{\mu_2}\cdot (T\varphi)_\mu (g_1z_1)\\
& = \tau_{\mu_2}^{-1}\cdot (T\varphi)_{\mu/G}\cdot \tau_{\mu_1}(g_1z_1)
 = \tau_{\mu_2}^{-1}\cdot (T\varphi)_{\mu/G}[z_1]= \tau_{\mu_2}^{-1}\cdot [z_2]\\
& = \tau_{\mu_2}^{-1}\cdot \tau_{\mu_2}(g_2z_2)= \Phi_{2g_2}(z_2)
= \Phi_{2g_2}\cdot (T\varphi)_\mu(z_1).
\end{align*}
Here we  denote by $\tau_{\mu_1}^{-1}(S)$ the pre-image of a subset
$S \subset (TQ_1)_{\mu_1}$ for the map
$\tau_{\mu_1}: (\mathbf{J}_L)_1^{-1}(\mu_1) \rightarrow
(TQ_1)_{\mu_1}$, and for any $z_1 \in
(\mathbf{J}_L)_1^{-1}(\mu_1)$, $\tau_{\mu_1}^{-1}\cdot \tau_{\mu_1}(z_1)=z_1$.
So, we obtain that $(T\varphi)_\mu \cdot
\Phi_{1g_1}=\Phi_{2g_2}\cdot (T\varphi)_\mu $. Moreover, we have
that
\begin{align*}
\mathcal{C}^L_2\cap (\mathbf{J}_L)_2^{-1}(\mu_2) & =
\tau_{\mu_2}^{-1}\cdot \tau_{\mu_2}(\mathcal{C}^L_2\cap (\mathbf{J}_L)_2^{-1}(\mu_2))
= \tau_{\mu_2}^{-1}\cdot \mathcal{C}^L_{2\mu_2}
= \tau_{\mu_2}^{-1}\cdot (T\varphi)_{\mu/G}(\mathcal{C}^L_{1\mu_1})\\
& = \tau_{\mu_2}^{-1}\cdot (T\varphi)_{\mu/G}\cdot
\tau_{\mu_1}(\mathcal{C}^L_1\cap (\mathbf{J}_L)_1^{-1}(\mu_1))
= \tau_{\mu_2}^{-1}\cdot \tau_{\mu_2}\cdot
(T\varphi)_\mu (\mathcal{C}^L_1\cap (\mathbf{J}_L)_1^{-1}(\mu_1))\\
& = (T\varphi)_\mu (\mathcal{C}^L_1\cap (\mathbf{J}_L)_1^{-1}(\mu_1)).
\end{align*}
Thus, the condition RpCL-1 holds. In the following we shall prove
that the condition RpCL-2 holds. For the above
$R_p$-reduced control laws $u^L_{i\mu_i}:(T
Q_i)_{\mu_i}\rightarrow \mathcal{C}^L_{i\mu_i}$, $i=1,2, $ there exist
the control laws $u^L_i: TQ_i \rightarrow \mathcal{C}^L_i, $  such that
$u^L_{i\mu_i}\cdot \tau_{\mu_i}=\tau_{\mu_i}\cdot u^L_i\cdot
j_{\mu_i}, \; i=1,2 $.
we shall prove that
$$ \xi_{(TQ_2,G_2,\omega^L_2, L_2,F^L_2, u^L_2)}\cdot T\varphi
=T(T\varphi)\cdot \xi_{(TQ_1,G_1,\omega^L_1,L_1,F^L_1,u^L_1)}.$$
In fact, from (4.4) we have that
\begin{align*}
& T((T\varphi)_{\mu/G})\cdot
\xi_{((TQ_1)_{\mu_1}, \omega^L_{1\mu_1}, l_{1\mu_1},
f^L_{1\mu_1},u^L_{1\mu_1})}\cdot \tau_{\mu_1}
= T((T\varphi)_{\mu/G})\cdot T\tau_{\mu_1}\cdot
\xi_{(TQ_1,G_1,\omega^L_1,L_1,F^L_1,u^L_1)}\cdot j_{\mu_1}\\
& =  T((T\varphi)_{\mu/G} \cdot \tau_{\mu_1})\cdot
\xi_{(TQ_1,G_1,\omega^L_1,L_1,F^L_1,u^L_1)}\cdot j_{\mu_1}
= T(\tau_{\mu_2}\cdot (T\varphi)_\mu)\cdot
\xi_{(TQ_1, G_1, \omega^L_1, L_1, F^L_1, u^L_1)}\cdot j_{\mu_1}\\
& = T\tau_{\mu_2}\cdot T(T\varphi) \cdot
\xi_{(TQ_1, G_1, \omega^L_1, L_1, F^L_1,u^L_1)}\cdot j_{\mu_1}.
\end{align*}
On the other hand,
\begin{align*}
& \xi_{((TQ_2)_{\mu_2}, \omega^L_{2\mu_2}, l_{2\mu_2},
f^L_{2\mu_2},u^L_{2\mu_2})}\cdot (T\varphi)_{\mu/G}\cdot \tau_{\mu_1}
= \xi_{((TQ_2)_{\mu_2}, \omega^L_{2\mu_2}, l_{2\mu_2},
f^L_{2\mu_2},u^L_{2\mu_2})}\cdot \tau_{\mu_2}\cdot (T\varphi)_\mu \\
& =T\tau_{\mu_2}\cdot \xi_{(TQ_2, G_2,
\omega^L_2, L_2, F^L_2, u^L_2)}\cdot j_{\mu_2}\cdot (T\varphi)_\mu
=T\tau_{\mu_2}\cdot \xi_{(TQ_2, G_2,
\omega^L_2, L_2, F^L_2, u^L_2)}\cdot T\varphi \cdot j_{\mu_1}.
\end{align*}
From (4.5) we have that
$$
T\tau_{\mu_2}\cdot \xi_{(TQ_2, G_2,
\omega^L_2, L_2, F^L_2, u^L_2)}\cdot T\varphi \cdot j_{\mu_1}
= T\tau_{\mu_2}\cdot T(T\varphi) \cdot
\xi_{(TQ_1,G_1,\omega^L_1,L_1,F^L_1,u^L_1)}\cdot j_{\mu_1}.
$$
Note that the map $ j_{\mu_1}: (\mathbf{J}_L)_1^{-1}(\mu_1)\rightarrow TQ_1 $ is injective,
and $T\tau_{\mu_2}: T(\mathbf{J}_L)_2^{-1}(\mu_2) \rightarrow
T(TQ_2)_{\mu_2}$ is surjective, hence, we have that
$$\xi_{(TQ_2,G_2,\omega^L_2,L_2,F^L_2,u^L_2)}\cdot T\varphi
=T(T\varphi)\cdot \xi_{(TQ_1,G_1,\omega^L_1,L_1,F^L_1,u^L_1)}.$$
It follows that the theorem holds.
\hskip 0.3cm $\blacksquare$\\

It is worthy of noting that, when the external force and control of a
regular point reducible RCL system $(TQ,G,\omega^L,L,F^L,\mathcal{C}^L)$ are both zeros,
that is, $F^L=0 $ and $\mathcal{C}^L=\emptyset$, in this case the RCL system
is just a regular point reducible Lagrangian system $(TQ,G,\omega^L,L)$.
Then the following theorem explains the relationship between the equivalence
for the regular point reducible Lagrangian systems with symmetries and the
equivalence for the associated $R_p$-reduced Lagrangian systems.

\begin{theo}
Two regular point reducible Lagrangian systems $(TQ_i, G_i, \omega^L_i,
L_i )$, $i=1,2,$ are equivalent if and
only if the associated $R_p$-reduced Lagrangian systems
$((TQ_i)_{\mu_i},\omega^L_{i\mu_i}, l_{i\mu_i} ),\; i=1,2,$ are equivalent.
\end{theo}

\noindent {\bf Proof:} If two regular point reducible Lagrangian systems $(TQ_i, G_i, \omega^L_i,
L_i )$, $i=1,2,$ are equivalent, then there exists a diffeomorphism
$\varphi: Q_1\rightarrow Q_2$ such that $T\varphi: TQ_1\rightarrow
TQ_2 $ is symplectic with respect to their Lagrangian symplectic forms
$\omega^L_i, \; i=1,2$, that is,
$\omega^L_1=(T\varphi)^\ast\cdot\omega^L_2$, and for $\mu_i \in
\mathfrak{g}^\ast_i,\; i=1,2$, $(T\varphi)_\mu=j_{\mu_2}^{-1}\cdot
T\varphi \cdot j_{\mu_1}: (\mathbf{J}_L)_1^{-1}(\mu_1)\rightarrow
(\mathbf{J}_L)_2^{-1}(\mu_2)$ is
$(G_{1\mu_1},G_{2\mu_2})$-equivariant. From the above
commutative Diagram-4,
we can define a map $(T\varphi)_{\mu/G}:(TQ_1)_{\mu_1}\rightarrow
(TQ_2)_{\mu_2}$, such that $(T\varphi)_{\mu/G} \cdot
\tau_{\mu_1}=\tau_{\mu_2}\cdot (T\varphi)_\mu$. Since
$(T\varphi)_\mu: (\mathbf{J}_L)_1^{-1}(\mu_1)\rightarrow
(\mathbf{J}_L)_2^{-1}(\mu_2)$ is $(G_{1\mu_1},G_{2\mu_2})$-equivariant,
then $(T\varphi)_{\mu/G}$ is well-defined.
In order to prove that the associated $R_p$-reduced Lagrangian systems
$((TQ_i)_{\mu_i},\omega^L_{i\mu_i}, l_{i\mu_i} ),\; i=1,2,$ are equivalent,
in the following we shall show that $(T\varphi)_{\mu/G}$ is symplectic
with respect to their $R_p$-reduced Lagrangian symplectic forms
$\omega^L_{i\mu_i}, \; i=1,2$, that is,
$(T\varphi)_{\mu/G}^\ast\omega^L_{2\mu_2}=\omega^L_{1\mu_1}$.
In fact, since $T\varphi: TQ_1\rightarrow TQ_2 $ is symplectic with
respect to their Lagrangian symplectic forms, the map $(T\varphi)^\ast:
\Omega^2(TQ_2)\rightarrow \Omega^2(TQ_1)$ satisfies $(T\varphi)^\ast
\omega^L_2=\omega^L_1$. From (2.3) we know that,
$j_{\mu_i}^\ast\omega^L_i=\tau_{\mu_i}^\ast\omega^L_{i\mu_i}, \;
i=1,2$, from the following commutative Diagram-5,
\[
\begin{CD}
\Omega^2(TQ_2) @ > j_{\mu_2}^\ast >>
\Omega^2((\mathbf{J}_L)_2^{-1}(\mu_2)) @ < \tau_{\mu_2}^\ast << \Omega^2((TQ_2)_{\mu_2})\\
@V(T\varphi)^\ast VV @V(T\varphi)_\mu^\ast VV @V(T\varphi)_{\mu/G}^\ast VV\\
\Omega^2(TQ_1) @ > j_{\mu_1}^\ast >>
\Omega^2((\mathbf{J}_L)_1^{-1}(\mu_1)) @< \tau_{\mu_1}^\ast <<
\Omega^2((TQ_1)_{\mu_1})
\end{CD}
\]
$$\mbox{Diagram-5}$$
we have that
\begin{align*}
\tau_{\mu_1}^\ast \cdot(T\varphi)_{\mu/G}^\ast\omega^L_{2\mu_2}&
=((T\varphi)_{\mu/G}\cdot \tau_{\mu_1})^\ast\omega^L_{2\mu_2}
=(\tau_{\mu_2}\cdot (T\varphi)_{\mu})^\ast\omega^L_{2\mu_2}\\ &
=(j_{\mu_2}^{-1}\cdot T\varphi \cdot j_{\mu_1})^\ast
\cdot \tau_{\mu_2}^\ast\omega^L_{2\mu_2}\\
&=j_{\mu_1}^\ast\cdot(T\varphi)^\ast\cdot(j_{\mu_2}^{-1})^\ast \cdot
j_{\mu_2}^\ast \omega^L_2\\ &= j_{\mu_1}^\ast\cdot
(T\varphi)^\ast\omega^L_2
=j_{\mu_1}^\ast\omega^L_1=\tau_{\mu_1}^\ast\omega^L_{1\mu_1}.
\end{align*}
Notice that $\tau_{\mu_1}$ is surjective, thus,
$(T\varphi)_{\mu/G}^\ast\omega^L_{2\mu_2}=\omega^L_{1\mu_1}$.\\

Conversely, assume that the $R_p$-reduced Lagrangian systems
$((TQ_i)_{\mu_i},\omega^L_{i\mu_i},l_{i\mu_i})$,
$i=1,2,$ are equivalent, then there exists a diffeomorphism
$(T\varphi)_{\mu/G}: (TQ_1)_{\mu_1}\rightarrow (TQ_2)_{\mu_2}$,
which is symplectic with respect to their $R_p$-reduced Lagrangian symplectic forms
$\omega^L_{i\mu_i}, \; i=1,2$.
From the above commutative Diagram-4, we can
define a map $(T\varphi)_\mu:
(\mathbf{J}_L)^{-1}_1(\mu_1)\rightarrow
(\mathbf{J}_L)^{-1}_2(\mu_2)$, such that $\tau_{\mu_2}\cdot
(T\varphi)_\mu= (T\varphi)_{\mu/G}\cdot \tau_{\mu_1}, $ and the map
$T\varphi: TQ_1 \rightarrow TQ_2 ,$ such that $T\varphi \cdot
j_{\mu_1}=j_{\mu_2}\cdot (T\varphi)_{\mu}, $,
as well as a diffeomorphism $\varphi: Q_1\rightarrow Q_2,
$, whose tangent map is just $T\varphi: TQ_1 \rightarrow TQ_2$. From
definition of $(T\varphi)_\mu$, we know that $(T\varphi)_\mu$ is
$(G_{1\mu_1},G_{2\mu_2})$-equivariant.
In the following we shall show that $T\varphi$ is
symplectic with respect to the Lagrangian symplectic forms $\omega^L_i, \;
i=1,2$, that is, $\omega^L_1=(T\varphi)^\ast\cdot\omega^L_2$.
Because $(T\varphi)_{\mu/G}: (TQ_1)_{\mu_1}\rightarrow
(TQ_2)_{\mu_2}$ is symplectic with respect to their $R_p$-reduced
Lagrangian symplectic forms, the map $((T\varphi)_{\mu/G})^\ast:
\Omega^2((TQ_2)_{\mu_2})\rightarrow \Omega^2((TQ_1)_{\mu_1})$,
satisfies $((T\varphi)_{\mu/G})^\ast\cdot
\omega^L_{2\mu_2}=\omega^L_{1\mu_1}$. From (2.3) we know that,
$j_{\mu_i}^\ast
\cdot\omega^L_i=\tau_{\mu_i}^\ast\cdot\omega^L_{i\mu_i}$, $i=1,2$,
from the commutative Diagram-5, we have that
\begin{align*}
j_{\mu_1}^\ast\cdot\omega^L_1
&=\tau_{\mu_1}^\ast\cdot\omega^L_{1\mu_1} =\tau_{\mu_1}^\ast
\cdot(T\varphi)_{\mu/G}^\ast\cdot\omega^L_{2\mu_2}
=((T\varphi)_{\mu/G}\cdot \tau_{\mu_1})^\ast\cdot\omega^L_{2\mu_2}\\
& =(\tau_{\mu_2}\cdot (T\varphi)_\mu)^\ast\cdot\omega^L_{2\mu_2}
=(j_{\mu_2}^{-1}\cdot T\varphi \cdot j_{\mu_1})^\ast\cdot
\tau_{\mu_2}^\ast \cdot\omega^L_{2\mu_2}\\ &
=j_{\mu_1}^\ast\cdot(T\varphi)^\ast \cdot(j^{-1}_{\mu_2})^\ast\cdot
j_{\mu_2}^\ast\cdot\omega^L_2=j_{\mu_1}^\ast\cdot
(T\varphi)^\ast\omega^L_2.
\end{align*}
Notice that $j_{\mu_1}$ is injective, and hence,
$\omega^L_1=(T\varphi)^\ast\omega^L_2$.
Thus, the regular point reducible Lagrangian systems $(TQ_i, G_i, \omega^L_i,
L_i )$, $i=1,2,$ are equivalent.
\hskip 1cm $\blacksquare $ \\

Thus, the regular point reduction Theorem 4.3 for the RCL systems
can be regarded as an extension of the regular
point reduction Theorem 4.4 for the regular Lagrangian systems under regular
controlled Lagrangian equivalence conditions. \\

\begin{rema}
If $(TQ, \omega^L)$ is a connected symplectic manifold, and
$\mathbf{J}_L: TQ \rightarrow \mathfrak{g}^\ast$ is a
non-equivariant momentum map with a non-equivariance group
one-cocycle $\sigma: G\rightarrow \mathfrak{g}^\ast$, which is
defined by $\sigma(g):=\mathbf{J}_L(g\cdot
z)-\operatorname{Ad}^\ast_{g^{-1}}\mathbf{J}_L(z)$, where $g\in G$
and $z\in TQ$. Then we know that $\sigma$ produces a new affine
action $\Theta: G\times \mathfrak{g}^\ast \rightarrow
\mathfrak{g}^\ast $ defined by
$\Theta(g,\mu):=\operatorname{Ad}^\ast_{g^{-1}}\mu + \sigma(g)$,
where $\mu \in \mathfrak{g}^\ast$, with respect to which the given
momentum map $\mathbf{J}_L$ is equivariant. Assume that $G$ acts
freely and properly on $TQ$, and $\tilde{G}_\mu$ denotes the
isotropy subgroup of $\mu \in \mathfrak{g}^\ast$ relative to this
affine action $\Theta$ and $\mu$ is a regular value of
$\mathbf{J}_L$. Then the quotient space
$(TQ)_\mu=\mathbf{J}_L^{-1}(\mu)/\tilde{G}_\mu$ is also a symplectic
manifold with the symplectic form $\omega^L_\mu$ uniquely
characterized by $(2.3)$. In this case, we can also define the
regular point reducible RCL system
$(TQ,G,\omega^L,L,F^L,\mathcal{C}^L)$ and RpCL-equivalence, and
prove the regular point reduction theorem for the RCL system by using
the above similar way.
\end{rema}

\section{Regular Orbit Reduction of the RCL System}

Since the set of regular Lagrangian systems with symmetries
on $TQ$ is a subset of the set of RCL systems with symmetries on $TQ$.
If we first admit the regular orbit reduction of a regular Lagrangian system
with symmetry, then we may study the regular orbit
reduction of an RCL system with symmetry,
as an extension of the regular orbit reduction of
a regular Lagrangian system under the regular controlled Lagrangian equivalence
conditions. In order to do these, in this section
we consider the RCL system with symmetry and momentum map,
and first give the regular orbit reducible RCL system
and the RoCL-equivalence, then prove the regular orbit reduction theorems
for the RCL system and regular Lagrangian system.\\

Note that, if an RCL system with symmetry and momentum map is
regular orbit reducible, then the associated regular
Lagrangian system must be regular orbit reducible.
Thus, from Definition 2.5 and Theorem 2.8,
if the Legendre transformation $\mathcal{F}L: TQ
\rightarrow T^{\ast}Q $ is $(\Phi^T,\; \Phi^{T*})$-equivariant,
then we can introduce
a kind of regular orbit reducible RCL systems as follows.
\begin{defi}
(Regular Orbit Reducible RCL System) A 6-tuple $(TQ, G,
\omega^L,L,F^L,\mathcal{C}^L)$, where the hyperregular Lagrangian
$L: TQ\rightarrow \mathbb{R}$, the fiber-preserving map $F^L:
TQ\rightarrow TQ$ and the fiber submanifold $\mathcal{C}^L$ of $TQ$
are all $G$-invariant, is called a {\bf regular orbit reducible RCL
system}, if the Legendre transformation $\mathcal{F}L: TQ
\rightarrow T^{\ast}Q $ is $(\Phi^T,\; \Phi^{T*})$-equivariant,
and there exists an orbit $\mathcal{O}_\mu, \;
\mu\in\mathfrak{g}^\ast$, where $\mu$ is a regular value of the
momentum map $\mathbf{J}_L$, such that the regular orbit reduced
system, that is, the 5-tuple
$((TQ)_{\mathcal{O}_\mu},\omega^L_{\mathcal{O}_\mu},l_{\mathcal{O}_\mu},f^L_{\mathcal{O}_\mu},
\mathcal{C}^L_{\mathcal{O}_\mu})$, where
$(TQ)_{\mathcal{O}_\mu}=\mathbf{J}_L^{-1}(\mathcal{O}_\mu)/G$,
$\tau_{\mathcal{O}_\mu}^\ast \omega^L_{\mathcal{O}_\mu}
=j_{\mathcal{O}_\mu}^\ast\omega^L-(\mathbf{J}_L)_{\mathcal{O}_\mu}^\ast
\omega^{L+}_{\mathcal{O}_\mu}$, $l_{\mathcal{O}_\mu}\cdot
\tau_{\mathcal{O}_\mu} =L\cdot j_{\mathcal{O}_\mu}$,
$F^L(\mathbf{J}_L^{-1}(\mathcal{O}_\mu))\subset
\mathbf{J}_L^{-1}(\mathcal{O}_\mu)$, $f^L_{\mathcal{O}_\mu}\cdot
\tau_{\mathcal{O}_\mu}=\tau_{\mathcal{O}_\mu}\cdot F^L\cdot
j_{\mathcal{O}_\mu}$, and $\mathcal{C}^L \cap
\mathbf{J}_L^{-1}(\mathcal{O}_\mu)\neq \emptyset $,
$\mathcal{C}^L_{\mathcal{O}_\mu}=\tau_{\mathcal{O}_\mu}(\mathcal{C}^L
\cap \mathbf{J}_L^{-1}(\mathcal{O}_\mu))$, is an RCL system, which is
simply written as $R_o$-reduced RCL system. Where
$((TQ)_{\mathcal{O}_\mu},\omega^L_{\mathcal{O}_\mu})$ is the
$R_o$-reduced space, the function
$l_{\mathcal{O}_\mu}:(TQ)_{\mathcal{O}_\mu}\rightarrow \mathbb{R}$
is called the $R_o$-reduced Lagrangian, the fiber-preserving map
$f^L_{\mathcal{O}_\mu}:(TQ)_{\mathcal{O}_\mu} \rightarrow (T
Q)_{\mathcal{O}_\mu}$ is called the $R_o$-reduced (external) force map,
$\mathcal{C}^L_{\mathcal{O}_\mu}$ is a fiber submanifold of $(T
Q)_{\mathcal{O}_\mu}$, and is called the $R_o$-reduced control subset.
\end{defi}

It is worthy of noting that for the regular orbit reducible RCL system
$(TQ,G,\omega^L,L,F^L,\mathcal{C}^L)$, the $G$-invariant external
force map $F^L: TQ \rightarrow TQ $ has to satisfy the conditions
$F^L(\mathbf{J}_L^{-1}(\mathcal{O}_\mu))\subset
\mathbf{J}_L^{-1}(\mathcal{O}_\mu), $ and
$f^L_{\mathcal{O}_\mu}\cdot
\tau_{\mathcal{O}_\mu}=\tau_{\mathcal{O}_\mu}\cdot F^L\cdot
j_{\mathcal{O}_\mu}$, such that we can define the $R_o$-reduced external
force map $f^L_{\mathcal{O}_\mu}:(TQ)_{\mathcal{O}_\mu} \rightarrow
(TQ)_{\mathcal{O}_\mu}. $ The condition $\mathcal{C}^L \cap
\mathbf{J}_L^{-1}(\mathcal{O}_\mu)\neq \emptyset $ in above
definition makes that the  $G$-invariant control subset
$\mathcal{C}^L \cap \mathbf{J}_L^{-1}(\mathcal{O}_\mu)$ can be
reduced and the $R_o$-reduced control subset is
$\mathcal{C}^L_{\mathcal{O}_\mu}=
\tau_{\mathcal{O}_\mu}(\mathcal{C}^L \cap
\mathbf{J}_L^{-1}(\mathcal{O}_\mu))$.\\

Assume that the dynamical vector field $\xi_{(TQ,G,\omega^L,L,F^L,u^L)}$ of a given
regular orbit reducible RCL system $(TQ, G,\omega^L,
L,F^L,\mathcal{C}^L)$ with a control law $u^L \in \mathcal{C}^L$ can be expressed by
\begin{equation}\xi_{(TQ,G,\omega^L,L,F^L,u^L)}
=\xi_L+\textnormal{vlift}(F^L)+\textnormal{vlift}(u^L).\label{5.1}
\end{equation}
Then, for the regular orbit reducible RCL system we can also
introduce the regular orbit reducible controlled Lagrangian
equivalence (RoCL-equivalence) as follows.
\begin{defi}
(RoCL-equivalence) Suppose that we have two regular orbit reducible
RCL systems $(TQ_i, G_i, \omega^L_i, L_i, F^L_i, \mathcal{C}^L_i)$,
$i=1,2$, we say them to be RoCL-equivalent, or simply,\\ $(TQ_1,
G_1, \omega^L_1, L_1, F^L_1,
\mathcal{C}^L_1)\stackrel{RoCL}{\sim}(TQ_2, G_2, \omega^L_2, L_2,
F^L_2, \mathcal{C}^L_2)$, if there exists a diffeomorphism
$\varphi:Q_1\rightarrow Q_2$ such that the following regular
orbit reducible controlled Lagrangian matching conditions hold:\\

\noindent {\bf RoCL-1:} For $\mathcal{O}_{\mu_i},\; \mu_i\in
\mathfrak{g}^\ast_i$, the regular reducible orbits of RCL systems
$(TQ_i, G_i, \omega^L_i, L_i, F^L_i, \mathcal{C}^L_i)$, $i=1,2$, the
map $(T\varphi)_{\mathcal{O}_\mu}=j_{\mathcal{O}_{\mu_2}}^{-1}\cdot
T\varphi\cdot j_{\mathcal{O}_{\mu_1}}:
(\mathbf{J}_L)_1^{-1}(\mathcal{O}_{\mu_1})\rightarrow
(\mathbf{J}_L)_2^{-1}(\mathcal{O}_{\mu_2})$ is
$(G_1,G_2)$-equivariant, $\mathcal{C}^L_2\cap
(\mathbf{J}_L)_2^{-1}(\mathcal{O}_{\mu_2})=(T\varphi)_{\mathcal{O}_\mu}
(\mathcal{C}^L_1\cap (\mathbf{J}_L)_1^{-1}(\mathcal{O}_{\mu_1}))$,
and $(\mathbf{J}_L)_{1\mathcal{O}_{\mu_1}}^\ast\cdot
\omega^{L+}_{1\mathcal{O}_{\mu_1}}=((T\varphi)_{\mathcal{O}_\mu})^\ast
\cdot(\mathbf{J}_L)_{2\mathcal{O}_{\mu_2}}^\ast\cdot\omega^{L+}_{2\mathcal{O}_{\mu_2}},$
where $\mu=(\mu_1, \mu_2)$, and denote by
$j_{\mathcal{O}_{\mu_2}}^{-1}(S)$ the pre-image of a subset
$S\subset TQ_2$ for the map $j_{\mathcal{O}_{\mu_2}}:
(\mathbf{J}_L)_2^{-1}(\mathcal{O}_{\mu_2})\rightarrow TQ_2$.

\noindent {\bf RoCL-2:}
For each control law $u^{L}_1:
TQ_1 \rightarrow \mathcal{C}^{L}_1, $ there exists the control law $u^{L}_2:
TQ_2 \rightarrow \mathcal{C}^{L}_2, $  such that the two
closed-loop dynamical systems produce the same dynamical vector fields, that is,
$\xi_{(TQ_2,G_2,\omega^{L}_2,L_2,F^{L}_2,u^{L}_2)}\cdot T\varphi
=T(T\varphi)\cdot \xi_{(TQ_1,G_1,\omega^{L}_1,L_1,F^{L}_1,u^{L}_1)}$.
\end{defi}

It is worthy of noting that for the regular orbit reducible RCL
system, the induced equivalent map $T\varphi$ not only keeps
the equivariance of $G$-action on their regular orbits, but also keeps
the restriction of the $(+)$-symplectic structure on the regular orbit to
$\mathbf{J}_L^{-1}(\mathcal{O}_\mu)$. If a feedback
control law $u^L: TQ \rightarrow \mathcal{C}^L $ is chosen,
and $u^L \in \mathcal{C}^L \cap \mathbf{J}^{-1}_L(\mathcal{O}_\mu),$
and $\mathcal{C}^L \cap \mathbf{J}^{-1}_L(\mathcal{O}_\mu)\neq \emptyset $,
then the $R_o$-reduced control law $u^L_{\mathcal{O}_\mu}:
(TQ)_{\mathcal{O}_\mu} \rightarrow \mathcal{C}^L_{\mathcal{O}_\mu}
=\tau_{\mathcal{O}_\mu}(\mathcal{C}^L\cap
\mathbf{J}^{-1}_L(\mathcal{O}_\mu))$, and
$u^L_{\mathcal{O}_\mu}\cdot\tau_{\mathcal{O}_\mu}
=\tau_{\mathcal{O}_\mu}\cdot u^L\cdot j_{\mathcal{O}_\mu}$.
The $R_o$-reduced RCL system $((TQ)_{\mathcal{O}_\mu},
\omega^L_{\mathcal{O}_\mu},l_{\mathcal{O}_\mu},
f^L_{\mathcal{O}_\mu}, u^L_{\mathcal{O}_\mu})$
is a closed-loop regular dynamical system with the $R_o$-reduced control law
$u^L_{\mathcal{O}_\mu}$. Assume that its dynamical vector field $\xi_{((T
Q)_{\mathcal{O}_\mu}, \omega^L_{\mathcal{O}_\mu},
l_{\mathcal{O}_\mu},f^L_{\mathcal{O}_\mu},u^L_{\mathcal{O}_\mu})}$
can be expressed by
\begin{equation}\xi_{((TQ)_{\mathcal{O}_\mu},
\omega^L_{\mathcal{O}_\mu},l_{\mathcal{O}_\mu},f^L_{\mathcal{O}_\mu},u^L_{\mathcal{O}_\mu})}=
\xi_{l_{\mathcal{O}_\mu}}+\textnormal{vlift}(f^L_{\mathcal{O}_\mu})
+\textnormal{vlift}(u^L_{\mathcal{O}_\mu}),
\label{5.2}\end{equation} where
$\xi_{l_{\mathcal{O}_\mu}}$ is the $R_o$-reduced Euler-Lagrange vector field,
$\textnormal{vlift}(f^L_{\mathcal{O}_\mu})=
\textnormal{vlift}(f^L_{\mathcal{O}_\mu})\xi_{l_{\mathcal{O}_\mu}}$,
$\textnormal{vlift}(u^L_{\mathcal{O}_\mu})=
\textnormal{vlift}(u^L_{\mathcal{O}_\mu})\xi_{l_{\mathcal{O}_\mu}}$
are the changes of $\xi_{l_{\mathcal{O}_\mu}}$ under the actions of the
$R_o$-reduced external force $f^L_{\mathcal{O}_\mu}$
and the $R_o$-reduced control law $u^L_{\mathcal{O}_\mu}$,
and the dynamical vector fields of the RCL system and the $R_o$-reduced
RCL system satisfy the condition
\begin{equation}\xi_{((TQ)_{\mathcal{O}_\mu},
\omega^L_{\mathcal{O}_\mu},l_{\mathcal{O}_\mu},f^L_{\mathcal{O}_\mu},u^L_{\mathcal{O}_\mu})}\cdot
\tau_{\mathcal{O}_\mu} =T\tau_{\mathcal{O}_\mu} \cdot \xi_{(T
Q,G,\omega^L,L,F^L,u^L)}\cdot
j_{\mathcal{O}_\mu},\label{5.3}
\end{equation}
see Marsden et al. \cite{mawazh10}, Wang \cite{wa18}.
Then we can obtain the following regular orbit reduction theorem for
the RCL system, which explains the relationship between the
RoCL-equivalence for the regular orbit reducible RCL system with
symmetry and the RCL-equivalence for the associated $R_o$-reduced RCL
system.
\begin{theo}
If two regular orbit reducible RCL systems $(TQ_i, G_i, \omega^L_i,
L_i, F^L_i, \mathcal{C}^L_i)$, $i=1,2,$ are RoCL-equivalent if and only if the
associated $R_o$-reduced RCL systems $((T
Q_i)_{\mathcal{O}_{\mu_i}}, \omega^L_{i\mathcal{O}_{\mu_i}},
l_{i\mathcal{O}_{\mu_i}}, f^L_{i\mathcal{O}_{\mu_i}}, \\
\mathcal{C}^L_{i\mathcal{O}_{\mu_i}})$, $i=1,2,$ are
RCL-equivalent.
\end{theo}

\noindent {\bf Proof:} If
$(TQ_1,G_1,\omega^L_1,L_1,F^L_1,\mathcal{C}^L_1)\stackrel{RoCL}{\sim}(T
Q_2,G_2,\omega^L_2,L_2,F^L_2,\mathcal{C}^L_2)$, then there exists a
diffeomorphism $\varphi: Q_1\rightarrow Q_2$, such that
for $\mathcal{O}_{\mu_i},\; \mu_i\in
\mathfrak{g}^\ast_i$, the regular reducible orbits, the
map $(T\varphi)_{\mathcal{O}_\mu}=j_{\mathcal{O}_{\mu_2}}^{-1}\cdot
T\varphi\cdot j_{\mathcal{O}_{\mu_1}}:
(\mathbf{J}_L)_1^{-1}(\mathcal{O}_{\mu_1})\rightarrow
(\mathbf{J}_L)_2^{-1}(\mathcal{O}_{\mu_2})$ is
$(G_1,G_2)$-equivariant, and $\mathcal{C}^L_2\cap
(\mathbf{J}_L)_2^{-1}(\mathcal{O}_{\mu_2})=(T\varphi)_{\mathcal{O}_\mu}
(\mathcal{C}^L_1\cap (\mathbf{J}_L)_1^{-1}(\mathcal{O}_{\mu_1}))$,
and RoCL-2 holds. From the following
commutative Diagram-6:
\[
\begin{CD}
TQ_1 @<j_{\mathcal{O}_{\mu_1}}<<( \mathbf{J}_L)_1^{-1}(\mathcal{O}_{\mu_1}) @>\tau_{\mathcal{O}_{\mu_1}}>> (TQ_1)_{\mathcal{O}_{\mu_1}}\\
@V T\varphi VV @V (T\varphi)_{\mathcal{O}_\mu} VV @V (T\varphi)_{\mathcal{O}_\mu/G}VV\\
TQ_2 @<j_{\mathcal{O}_{\mu_2}}<< (\mathbf{J}_L)_2^{-1}(\mathcal{O}_{\mu_2})
@>\tau_{\mathcal{O}_{\mu_2}}>>(TQ_2)_{\mathcal{O}_{\mu_2}}
\end{CD}
\]
$$\mbox{Diagram-6}$$
we can define a map $(T\varphi)_{\mathcal{O}_\mu/G}:(T
Q_1)_{\mathcal{O}_{\mu_1}}\rightarrow (TQ_2)_{\mathcal{O}_{\mu_2}}$ such that
$(T\varphi)_{\mathcal{O}_\mu/G}\cdot
\tau_{\mathcal{O}_{\mu_1}}=\tau_{\mathcal{O}_{\mu_2}}\cdot (T\varphi)_{\mathcal{O}_\mu}$.
Because $(T\varphi)_{\mathcal{O}_\mu}: (\mathbf{J}_L)_1^{-1}(\mathcal{O}_{\mu_1})\rightarrow
(\mathbf{J}_L)_2^{-1}(\mathcal{O}_{\mu_2})$ is $(G_{1},G_{2})$-equivariant,
$(T\varphi)_{\mathcal{O}_\mu/G}$ is well-defined. We shall show that
$\mathcal{C}^L_{2\mathcal{O}_{\mu_2}}
=(T\varphi)_{\mathcal{O}_\mu/G}(\mathcal{C}^L_{1\mathcal{O}_{\mu_1}})$. In fact,
because $(TQ_i, G_i, \omega^L_i,
L_i, F^L_i,\mathcal{C}^L_i)$, $i=1,2,$ are regular orbit reducible RCL systems,
then $\mathcal{C}^L_i\cap (\mathbf{J}_L)_i^{-1}(\mathcal{O}_{\mu_i})\neq
\emptyset $ and $\mathcal{C}^L_{i\mathcal{O}_{\mu_i}}=\tau_{\mathcal{O}_{\mu_i}}(\mathcal{C}^L_i\cap
(\mathbf{J}_L)_i^{-1}(\mathcal{O}_{\mu_i})),\; i=1,2$. From $\mathcal{C}^L_2\cap
(\mathbf{J}_L)_2^{-1}(\mathcal{O}_{\mu_2})=(T\varphi)_{\mathcal{O}_\mu} (\mathcal{C}^L_1\cap
(\mathbf{J}_L)_1^{-1}(\mathcal{O}_{\mu_1}))$, we have that
\begin{align*}
\mathcal{C}^L_{2\mathcal{O}_{\mu_2}} & =\tau_{\mathcal{O}_{\mu_2}}(\mathcal{C}^L_2\cap (\mathbf{J}_L)_2^{-1}(\mathcal{O}_{\mu_2}))
=\tau_{\mathcal{O}_{\mu_2}}\cdot (T\varphi)_{\mathcal{O}_\mu} (\mathcal{C}^L_1\cap
(\mathbf{J}_L)_1^{-1}(\mathcal{O}_{\mu_1}))\\ & =(T\varphi)_{\mathcal{O}_\mu/G}\cdot
\tau_{\mathcal{O}_{\mu_1}}(\mathcal{C}^L_1\cap (\mathbf{J}_L)_1^{-1}(\mathcal{O}_{\mu_1}))
=(T\varphi)_{\mathcal{O}_\mu/G}(\mathcal{C}^L_{1\mathcal{O}_{\mu_1}}).
\end{align*}
Thus, the condition RCL-1 holds. On the other hands,
for the $R_o$-reduced control law $u^L_{1\mathcal{O}_{\mu_1}}:(T
Q_1)_{\mathcal{O}_{\mu_1}}\rightarrow \mathcal{C}^L_{1\mathcal{O}_{\mu_1}}$, we have the control law
$u^L_1: TQ_1 \rightarrow \mathcal{C}^L_1, $, such that
$u^L_{1\mathcal{O}_{\mu_1}}\cdot\tau_{\mathcal{O}_{\mu_1}}=\tau_{\mathcal{O}_{\mu_1}}\cdot u^L_1\cdot
j_{\mathcal{O}_{\mu_1}}$. From the condition RoCL-2 we know that
there exists the control law $u^L_2:
TQ_2 \rightarrow \mathcal{C}^L_2, $  such that
$\xi_{(TQ_2,G_2,\omega^L_2,L_2,F^L_2,u^L_2)}\cdot T\varphi
=T(T\varphi)\cdot \xi_{(TQ_1,G_1,\omega^L_1,L_1,F^L_1,u^L_1)}$.
But, for the control law $u^L_2: TQ_2 \rightarrow \mathcal{C}^L_2, $ we have
the $R_o$-reduced control law $u^L_{2\mathcal{O}_{\mu_2}}:(T
Q_2)_{\mathcal{O}_{\mu_2}}\rightarrow \mathcal{C}^L_{2\mathcal{O}_{\mu_2}}$, such that
$u^L_{2\mathcal{O}_{\mu_2}}\cdot\tau_{\mathcal{O}_{\mu_2}}=\tau_{\mathcal{O}_{\mu_2}}\cdot u^L_2\cdot
j_{\mathcal{O}_{\mu_2}}$. Note that for $i=1,2$,
from (5.3), we have that
\begin{align}
\xi_{((TQ_i)_{\mathcal{O}_{\mu_i}}, \omega^L_{i\mathcal{O}_{\mu_i}}, l_{i\mathcal{O}_{\mu_i}},
f^L_{i\mathcal{O}_{\mu_i}},u^L_{i\mathcal{O}_{\mu_i}})}\cdot \tau_{\mathcal{O}_{\mu_i}}=T\tau_{\mathcal{O}_{\mu_i}}\cdot
\xi_{(TQ_i, G_i, \omega^L_i, L_i, F^L_i, u^L_i)}\cdot j_{\mathcal{O}_{\mu_i}},
\end{align}
and from the commutative Diagram-6,
$(T\varphi)_{\mathcal{O}_\mu/G}\cdot
\tau_{\mathcal{O}_{\mu_1}}=\tau_{\mathcal{O}_{\mu_2}}\cdot (T\varphi)_{\mathcal{O}_\mu}$ and
$j_{\mathcal{O}_{\mu_2}}\cdot (T\varphi)_{\mathcal{O}_\mu} = (T\varphi) \cdot j_{\mathcal{O}_{\mu_1}}$,
then we have that
\begin{align*}
& \xi_{((TQ_2)_{\mathcal{O}_{\mu_2}}, \omega^L_{2\mathcal{O}_{\mu_2}}, l_{2\mathcal{O}_{\mu_2}},
f^L_{2\mathcal{O}_{\mu_2}},u^L_{2\mathcal{O}_{\mu_2}})}\cdot (T\varphi)_{\mathcal{O}_\mu/G}\cdot \tau_{\mathcal{O}_{\mu_1}}\\
& = \xi_{((TQ_2)_{\mathcal{O}_{\mu_2}}, \omega^L_{2\mathcal{O}_{\mu_2}}, l_{2\mathcal{O}_{\mu_2}},
f^L_{2\mathcal{O}_{\mu_2}},u^L_{2\mathcal{O}_{\mu_2}})}\cdot \tau_{\mathcal{O}_{\mu_2}}\cdot (T\varphi)_{\mathcal{O}_\mu} \\
& =T\tau_{\mathcal{O}_{\mu_2}}\cdot \xi_{(TQ_2, G_2,
\omega^L_2, L_2, F^L_2, u^L_2)}\cdot j_{\mathcal{O}_{\mu_2}}\cdot (T\varphi)_{\mathcal{O}_\mu}\\
& =T\tau_{\mathcal{O}_{\mu_2}}\cdot \xi_{(TQ_2, G_2,
\omega^L_2, L_2, F^L_2, u^L_2)}\cdot (T\varphi) \cdot j_{\mathcal{O}_{\mu_1}}\\
& = T\tau_{\mathcal{O}_{\mu_2}}\cdot T(T\varphi) \cdot
\xi_{(TQ_1,G_1,\omega^L_1,L_1,F^L_1,u^L_1)}\cdot j_{\mathcal{O}_{\mu_1}}\\
& = T(\tau_{\mathcal{O}_{\mu_2}}\cdot (T\varphi)_{\mathcal{O}_\mu})\cdot
\xi_{(TQ_1,G_1,\omega^L_1,L_1,F^L_1,u^L_1)}\cdot j_{\mathcal{O}_{\mu_1}}\\
& = T((T\varphi)_{\mathcal{O}_\mu/G}\cdot \tau_{\mathcal{O}_{\mu_1}})\cdot
\xi_{(TQ_1,G_1,\omega^L_1,L_1,F^L_1,u^L_1)}\cdot j_{\mathcal{O}_{\mu_1}}\\
& =T((T\varphi)_{\mathcal{O}_\mu/G})\cdot T\tau_{\mathcal{O}_{\mu_1}}\cdot
\xi_{(TQ_1,G_1,\omega^L_1,L_1,F^L_1,u^L_1)}\cdot j_{\mathcal{O}_{\mu_1}}\\
& = T((T\varphi)_{\mathcal{O}_\mu/G})\cdot
\xi_{((TQ_1)_{\mathcal{O}_{\mu_1}}, \omega^L_{1\mathcal{O}_{\mu_1}}, l_{1\mathcal{O}_{\mu_1}},
f^L_{1\mathcal{O}_{\mu_1}},u^L_{1\mathcal{O}_{\mu_1}})}\cdot \tau_{\mathcal{O}_{\mu_1}}.
\end{align*}
Since $\tau_{\mathcal{O}_{\mu_1}}: (\mathbf{J}_L)_1^{-1}(\mathcal{O}_{\mu_1}) \rightarrow
(TQ_1)_{\mathcal{O}_{\mu_1}}$ is surjective, thus,
\begin{align}
& \xi_{((TQ_2)_{\mathcal{O}_{\mu_2}}, \omega^L_{2\mathcal{O}_{\mu_2}}, l_{2\mathcal{O}_{\mu_2}},
f^L_{2\mathcal{O}_{\mu_2}},u^L_{2\mathcal{O}_{\mu_2}})}\cdot (T\varphi)_{\mathcal{O}_\mu/G}\nonumber \\
& =T((T\varphi)_{\mathcal{O}_\mu/G})\cdot
\xi_{((TQ_1)_{\mathcal{O}_{\mu_1}}, \omega^L_{1\mathcal{O}_{\mu_1}}, l_{1\mathcal{O}_{\mu_1}},
f^L_{1\mathcal{O}_{\mu_1}},u^L_{1\mathcal{O}_{\mu_1}})},
\end{align}
that is, the condition RCL-2 holds.
So, the $R_o$-reduced RCL systems $((T
Q_i)_{\mathcal{O}_{\mu_i}},\omega^L_{i\mathcal{O}_{\mu_i}},l_{i\mathcal{O}_{\mu_i}}, f^L_{i\mathcal{O}_{\mu_i}}, \\
\mathcal{C}^L_{i\mathcal{O}_{\mu_i}}),$ $i=1,2,$ are RCL-equivalent.\\

Conversely, assume that the $R_o$-reduced RCL systems $((T
Q_i)_{\mathcal{O}_{\mu_i}},\omega^L_{i\mathcal{O}_{\mu_i}},l_{i\mathcal{O}_{\mu_i}},
f^L_{i\mathcal{O}_{\mu_i}},\mathcal{C}^L_{i\mathcal{O}_{\mu_i}})$,
$i=1,2,$ are RCL-equivalent, then there exists a diffeomorphism
$(T\varphi)_{\mathcal{O}_\mu/G}: (TQ_1)_{\mathcal{O}_{\mu_1}}\rightarrow (T
Q_2)_{\mathcal{O}_{\mu_2}}$, such that
$\mathcal{C}^L_{2\mathcal{O}_{\mu_2}}=(T\varphi)_{\mathcal{O}_\mu/G}(\mathcal{C}^L_{1\mathcal{O}_{\mu_1}}),\;
\forall \mathcal{O}_{\mu_i}, \; \mu_i\in\mathfrak{g}_i^\ast, \; i=1,2$ and
for the $R_o$-reduced control law $u^L_{1\mathcal{O}_{\mu_1}}:(T
Q_1)_{\mathcal{O}_{\mu_1}}\rightarrow \mathcal{C}^L_{1\mathcal{O}_{\mu_1}}$, there exists
the $R_o$-reduced control law $u^L_{2\mathcal{O}_{\mu_2}}:(T
Q_2)_{\mathcal{O}_{\mu_2}}\rightarrow \mathcal{C}^L_{2\mathcal{O}_{\mu_2}}$, such that (5.5) holds.
Then from commutative Diagram-6, we can
define a map $(T\varphi)_{\mathcal{O}_\mu}: (\mathbf{J}_L)^{-1}_1(\mathcal{O}_{\mu_1})\rightarrow
(\mathbf{J}_L)^{-1}_2(\mathcal{O}_{\mu_2})$ such that $\tau_{\mathcal{O}_{\mu_2}}\cdot
(T\varphi)_{\mathcal{O}_\mu}=(T\varphi)_{\mathcal{O}_\mu/G}\cdot \tau_{\mathcal{O}_{\mu_1}}, $ and the
map $T\varphi: TQ_1\rightarrow TQ_2 $ such that
$T\varphi \cdot j_{\mathcal{O}_{\mu_1}}=j_{\mathcal{O}_{\mu_2}}\cdot (T\varphi)_{\mathcal{O}_{\mu}}, $
as well as a diffeomorphism $\varphi:
Q_1\rightarrow Q_2, $ whose tangent lift is just
$ T\varphi: TQ_1 \rightarrow TQ_2$. Moreover, for above definition of
$(T\varphi)_{\mathcal{O}_\mu}$, we know that $(T\varphi)_{\mathcal{O}_\mu}$ is
$(G_{1},G_{2})$-equivariant. In fact, for any $z_i\in
(\mathbf{J}_L)_i^{-1}(\mathcal{O}_{\mu_i})$, $g_i \in G_{i}$, $i=1,2$ such that
$z_2=(T\varphi)_{\mathcal{O}_\mu} (z_1)$, $[z_2]=(T\varphi)_{\mathcal{O}_\mu/G}[z_1]$, then
we have that
\begin{align*}
(T\varphi)_{\mathcal{O}_\mu }(\Phi_{1g_1}(z_1))& = \tau_{\mathcal{O}_{\mu_2}}^{-1}\cdot
\tau_{\mathcal{O}_{\mu_2}}\cdot (T\varphi)_{\mathcal{O}_\mu} (\Phi_{1g_1}(z_1))
= \tau_{\mathcal{O}_{\mu_2}}^{-1}\cdot \tau_{\mathcal{O}_{\mu_2}}\cdot (T\varphi)_{\mathcal{O}_\mu} (g_1z_1)\\
& = \tau_{\mathcal{O}_{\mu_2}}^{-1}\cdot (T\varphi)_{\mathcal{O}_\mu/G}\cdot \tau_{\mathcal{O}_{\mu_1}}(g_1z_1)
 = \tau_{\mathcal{O}_{\mu_2}}^{-1}\cdot (T\varphi)_{\mathcal{O}_\mu/G}[z_1]= \tau_{\mathcal{O}_{\mu_2}}^{-1}\cdot [z_2]\\
& = \tau_{\mathcal{O}_{\mu_2}}^{-1}\cdot \tau_{\mathcal{O}_{\mu_2}}(g_2z_2)= \Phi_{2g_2}(z_2)
= \Phi_{2g_2}\cdot (T\varphi)_{\mathcal{O}_\mu}(z_1).
\end{align*}
Here we  denote by $\tau_{\mathcal{O}_{\mu_1}}^{-1}(S)$ the pre-image of a subset
$S \subset (TQ_1)_{\mathcal{O}_{\mu_1}}$ for the map
$\tau_{\mathcal{O}_{\mu_1}}: (\mathbf{J}_L)_1^{-1}(\mathcal{O}_{\mu_1}) \rightarrow
(TQ_1)_{\mathcal{O}_{\mu_1}}$, and for any $z_1 \in
(\mathbf{J}_L)_1^{-1}(\mathcal{O}_{\mu_1})$, $\tau_{\mathcal{O}_{\mu_1}}^{-1}\cdot \tau_{\mathcal{O}_{\mu_1}}(z_1)=z_1$.
So, we obtain that $(T\varphi)_{\mathcal{O}_\mu} \cdot
\Phi_{1g_1}=\Phi_{2g_2}\cdot (T\varphi)_{\mathcal{O}_\mu} $.
Moreover, we have that
\begin{align*}
\mathcal{C}^L_2\cap (\mathbf{J}_L)_2^{-1}(\mathcal{O}_{\mu_2}) & =
\tau_{\mathcal{O}_{\mu_2}}^{-1}\cdot \tau_{\mathcal{O}_{\mu_2}}(\mathcal{C}^L_2\cap (\mathbf{J}_L)_2^{-1}(\mathcal{O}_{\mu_2})) \\
& = \tau_{\mathcal{O}_{\mu_2}}^{-1}\cdot \mathcal{C}^L_{2\mathcal{O}_{\mu_2}}
= \tau_{\mathcal{O}_{\mu_2}}^{-1}\cdot (T\varphi)_{\mathcal{O}_\mu/G}(\mathcal{C}^L_{1\mathcal{O}_{\mu_1}})\\
& = \tau_{\mathcal{O}_{\mu_2}}^{-1}\cdot (T\varphi)_{\mathcal{O}_\mu/G}\cdot
\tau_{\mathcal{O}_{\mu_1}}(\mathcal{C}^L_1\cap (\mathbf{J}_L)_1^{-1}(\mathcal{O}_{\mu_1})) \\
& = \tau_{\mathcal{O}_{\mu_2}}^{-1}\cdot \tau_{\mathcal{O}_{\mu_2}}\cdot
(T\varphi)_{\mathcal{O}_\mu }(\mathcal{C}^L_1\cap (\mathbf{J}_L)_1^{-1}(\mathcal{O}_{\mu_1}))\\
& = (T\varphi)_{\mathcal{O}_\mu} (\mathcal{C}^L_1\cap (\mathbf{J}_L)_1^{-1}(\mathcal{O}_{\mu_1})).
\end{align*}
Thus, the condition RoCL-1 holds. In the following we shall prove
that the condition RoCL-2 holds. For the above
$R_o$-reduced control laws $u^L_{i\mathcal{O}_{\mu_i}}:(T
Q_i)_{\mathcal{O}_{\mu_i}}\rightarrow \mathcal{C}^L_{i\mathcal{O}_{\mu_i}}$, $i=1,2, $ there exist
the control laws $u^L_i: TQ_i \rightarrow \mathcal{C}^L_i, $  such that
$u^L_{i\mathcal{O}_{\mu_i}}\cdot \tau_{\mathcal{O}_{\mu_i}}=\tau_{\mathcal{O}_{\mu_i}}\cdot u^L_i\cdot
j_{\mathcal{O}_{\mu_i}}, \; i=1,2 $.
we shall prove that
$$ \xi_{(TQ_2,G_2,\omega^L_2, L_2,F^L_2, u^L_2)}\cdot T\varphi
=T(T\varphi)\cdot \xi_{(TQ_1,G_1,\omega^L_1,L_1,F^L_1,u^L_1)}.$$
In fact, from (5.4) we have that
\begin{align*}
& T((T\varphi)_{\mathcal{O}_\mu/G})\cdot
\xi_{((TQ_1)_{\mathcal{O}_{\mu_1}}, \omega^L_{1\mathcal{O}_{\mu_1}}, l_{1\mathcal{O}_{\mu_1}},
f^L_{1\mathcal{O}_{\mu_1}},u^L_{1\mathcal{O}_{\mu_1}})}\cdot \tau_{\mathcal{O}_{\mu_1}}\\
& = T((T\varphi)_{\mathcal{O}_\mu/G})\cdot T\tau_{\mathcal{O}_{\mu_1}}\cdot
\xi_{(TQ_1,G_1,\omega^L_1,L_1,F^L_1,u^L_1)}\cdot j_{\mathcal{O}_{\mu_1}}\\
& =  T((T\varphi)_{\mathcal{O}_\mu/G} \cdot \tau_{\mathcal{O}_{\mu_1}})\cdot
\xi_{(TQ_1,G_1,\omega^L_1,L_1,F^L_1,u^L_1)}\cdot j_{\mathcal{O}_{\mu_1}}\\
& = T(\tau_{\mathcal{O}_{\mu_2}}\cdot (T\varphi)_{\mathcal{O}_\mu})\cdot
\xi_{(TQ_1, G_1, \omega^L_1, L_1, F^L_1, u^L_1)}\cdot j_{\mathcal{O}_{\mu_1}}\\
& = T\tau_{\mathcal{O}_{\mu_2}}\cdot T(T\varphi) \cdot
\xi_{(TQ_1, G_1, \omega^L_1, L_1, F^L_1,u^L_1)}\cdot j_{\mathcal{O}_{\mu_1}}.
\end{align*}
On the other hand,
\begin{align*}
& \xi_{((TQ_2)_{\mathcal{O}_{\mu_2}}, \omega^L_{2\mathcal{O}_{\mu_2}}, l_{2\mathcal{O}_{\mu_2}},
f^L_{2\mathcal{O}_{\mu_2}},u^L_{2\mathcal{O}_{\mu_2}})}\cdot (T\varphi)_{\mathcal{O}_\mu/G}\cdot \tau_{\mathcal{O}_{\mu_1}}\\
& = \xi_{((TQ_2)_{\mathcal{O}_{\mu_2}}, \omega^L_{2\mathcal{O}_{\mu_2}}, l_{2\mathcal{O}_{\mu_2}},
f^L_{2\mathcal{O}_{\mu_2}},u^L_{2\mathcal{O}_{\mu_2}})}\cdot \tau_{\mathcal{O}_{\mu_2}}\cdot (T\varphi)_{\mathcal{O}_\mu} \\
& =T\tau_{\mathcal{O}_{\mu_2}}\cdot \xi_{(TQ_2, G_2,
\omega^L_2, L_2, F^L_2, u^L_2)}\cdot j_{\mathcal{O}_{\mu_2}}\cdot (T\varphi)_{\mathcal{O}_\mu}\\
& =T\tau_{\mathcal{O}_{\mu_2}}\cdot \xi_{(TQ_2, G_2,
\omega^L_2, L_2, F^L_2, u^L_2)}\cdot T\varphi \cdot j_{\mathcal{O}_{\mu_1}}.
\end{align*}
From (5.5) we have that
$$
T\tau_{\mathcal{O}_{\mu_2}}\cdot \xi_{(TQ_2, G_2,
\omega^L_2, L_2, F^L_2, u^L_2)}\cdot T\varphi \cdot j_{\mathcal{O}_{\mu_1}}
= T\tau_{\mathcal{O}_{\mu_2}}\cdot T(T\varphi) \cdot
\xi_{(TQ_1,G_1,\omega^L_1,L_1,F^L_1,u^L_1)}\cdot j_{\mathcal{O}_{\mu_1}}.
$$
Note that the map $ j_{\mathcal{O}_{\mu_1}}: (\mathbf{J}_L)_1^{-1}(\mathcal{O}_{\mu_1})\rightarrow TQ_1 $ is injective,
and $T\tau_{\mathcal{O}_{\mu_2}}: T(\mathbf{J}_L)_2^{-1}(\mathcal{O}_{\mu_2}) \rightarrow
T(TQ_2)_{\mathcal{O}_{\mu_2}}$ is surjective, hence, we have that
$$ \xi_{(TQ_2,G_2,\omega^L_2,L_2,F^L_2,u^L_2)}\cdot T\varphi
=T(T\varphi)\cdot \xi_{(TQ_1,G_1,\omega^L_1,L_1,F^L_1,u^L_1)}.$$
It follows that the theorem holds.
\hskip 0.3cm $\blacksquare$\\

It is worthy of noting that, when the external force and control of a
regular orbit reducible RCL system $(TQ,G,\omega^L,L,F^L,\mathcal{C}^L)$ are both zeros,
that is, $F^L=0 $ and $\mathcal{C}^L=\emptyset$, in this case the RCL system
is just a regular orbit reducible Lagrangian system $(TQ,G,\omega^L,L)$.
Then the following theorem explains the relationship between the equivalence
for the regular orbit reducible Lagrangian systems with symmetries and the
equivalence for the associated $R_o$-reduced Lagrangian systems.

\begin{theo}
If two regular orbit reducible Lagrangian systems $(TQ_i, G_i, \omega^L_i,
L_i )$, $i=1,2,$ are equivalent, then
their associated $R_o$-reduced Lagrangian systems $((T
Q)_{\mathcal{O}_{\mu_i}}, \omega^L_{i\mathcal{O}_{\mu_i}},
l_{i\mathcal{O}_{\mu_i}})$, $i=1,2,$ must be
equivalent. Conversely, if the $R_o$-reduced Lagrangian systems $((T
Q)_{\mathcal{O}_{\mu_i}}, \omega^L_{i\mathcal{O}_{\mu_i}},
l_{i\mathcal{O}_{\mu_i}})$, $i=1,2,$ are equivalent,
and the induced map
$(T\varphi)_{\mathcal{O}_\mu}:(\mathbf{J}_L)_1^{-1}(\mathcal{O}_{\mu_1})\rightarrow
(\mathbf{J}_L)_2^{-1}(\mathcal{O}_{\mu_2})$, such that
$(\mathbf{J}_L)_{1\mathcal{O}_{\mu_1}}^\ast\cdot
\omega^{L+}_{1\mathcal{O}_{\mu_1}}=(T\varphi)_{\mathcal{O}_\mu}^\ast
\cdot(\mathbf{J}_L)_{2\mathcal{O}_{\mu_2}}^\ast\cdot\omega^{L+}_{2\mathcal{O}_{\mu_2}},$
then the regular orbit reducible Lagrangian systems $(TQ_i, G_i,
\omega^L_i, L_i )$, $i=1,2,$ are equivalent.
\end{theo}

\noindent {\bf Proof:} If two regular orbit reducible Lagrangian
systems $(TQ_i, G_i, \omega^L_i, L_i )$, $i=1,2,$ are equivalent,
then there exists a diffeomorphism $\varphi: Q_1\rightarrow Q_2$,
such that $T\varphi: TQ_1\rightarrow TQ_2 $ is symplectic with respect to their
Lagrangian symplectic forms $\omega^L_i, \; i=1,2, $
and for $\mathcal{O}_{\mu_i}, \; \mu_i\in \mathfrak{g}_i^\ast, \; i=1,2 $,
$(T\varphi)_{\mathcal{O}_\mu}=j_{\mathcal{O}_{\mu_2}}^{-1}\cdot
T\varphi \cdot j_{\mathcal{O}_{\mu_1}}:
(\mathbf{J}_L)_1^{-1}(\mathcal{O}_{\mu_1})\rightarrow
(\mathbf{J}_L)_2^{-1}(\mathcal{O}_{\mu_2})$ is
$(G_1,G_2)$-equivariant. From the above
commutative Diagram-6,
we can define a map $(T\varphi)_{\mathcal{O}_\mu/G}:(T
Q_1)_{\mathcal{O}_{\mu_1}}\rightarrow (T
Q_2)_{\mathcal{O}_{\mu_2}}$, such that
$(T\varphi)_{\mathcal{O}_\mu/G} \cdot
\tau_{\mathcal{O}_{\mu_1}}=\tau_{\mathcal{O}_{\mu_2}}\cdot
(T\varphi)_{\mathcal{O}_\mu}$. Since
$(T\varphi)_{\mathcal{O}_\mu}:
(\mathbf{J}_L)_1^{-1}(\mathcal{O}_{\mu_1})\rightarrow
(\mathbf{J}_L)_2^{-1}(\mathcal{O}_{\mu_2})$ is
$(G_1,G_2)$-equivariant, then $(T\varphi)_{\mathcal{O}_\mu/G}$ is
well-defined. In order to prove that
the associated $R_o$-reduced Lagrangian systems $((T
Q)_{\mathcal{O}_{\mu_i}}, \omega^L_{i\mathcal{O}_{\mu_i}},
l_{i\mathcal{O}_{\mu_i}})$, $i=1,2,$ are equivalent,
in following we shall prove that $(T\varphi)_{\mathcal{O}_\mu/G}$ is
symplectic with respect to their $R_o$-reduced Lagrangian symplectic forms
$\omega^L_{i\mathcal{O}_i}, \; i=1,2$, that is,
$(T\varphi)_{\mathcal{O}_\mu/G}^\ast\cdot
\omega^L_{2\mathcal{O}_{\mu_2}}=\omega^L_{1\mathcal{O}_{\mu_1}}$.
In fact, since $T\varphi:
TQ_1\to TQ_2$ is symplectic with respect to their Lagrangian symplectic forms, and
the map $(T\varphi)^\ast: \Omega^2(TQ_2)\rightarrow \Omega^2(TQ_1)$
satisfies $(T\varphi)^\ast \cdot\omega^L_2=\omega^L_1$. From (2.4) we have that
$j_{\mathcal{O}_{\mu_i}}^\ast\cdot
\omega^L_i=\tau_{\mathcal{O}_{\mu_i}}^\ast\cdot\omega^L_{i\mathcal{O}_{\mu_i}}
+(\mathbf{J}_L)_{i\mathcal{O}_{\mu_i}}^\ast\cdot
\omega^{L+}_{i\mathcal{O}_{\mu_i}}$, $i=1,2,$ and
$(\mathbf{J}_L)_{1\mathcal{O}_{\mu_1}}^\ast\cdot\omega^{L+}_{1\mathcal{O}_{\mu_1}}
=((T\varphi)_{\mathcal{O}_\mu})^\ast
\cdot (\mathbf{J}_L)_{2\mathcal{O}_{\mu_2}}^\ast\cdot
\omega^{L+}_{2\mathcal{O}_{\mu_2}}$, from the following commutative
Diagram-7,
\[
\begin{CD}
\Omega^2(TQ_2) @> j_{\mathcal{O}_{\mu_2}}^\ast >>
\Omega^2((\mathbf{J}_L)_2^{-1}(\mathcal{O}_{\mu_2}))
@< \tau_{\mathcal{O}_{\mu_2}}^\ast << \Omega^2((TQ_2)_{\mathcal{O}_{\mu_2}})\\
@V (T\varphi)^\ast VV @ V(T\varphi)_{\mathcal{O}_\mu}^\ast
VV @V (T\varphi)_{\mathcal{O}_\mu/G}^\ast VV\\
\Omega^2(TQ_1) @> j_{\mathcal{O}_{\mu_1}}^\ast >>
\Omega^2((\mathbf{J}_L)_1^{-1}(\mathcal{O}_{\mu_1})) @<
\tau_{\mathcal{O}_{\mu_1}}^\ast <<
\Omega^2((TQ_1)_{\mathcal{O}_{\mu_1}})
\end{CD}\]
$$\mbox{Diagram-7}$$
we have that
\begin{align*}
\tau_{\mathcal{O}_{\mu_1}}^\ast\cdot
(T\varphi)_{\mathcal{O}_\mu/G}^\ast \omega^L_{2\mathcal{O}_{\mu_2}}
& =((T\varphi)_{\mathcal{O}_\mu/G}\cdot
\tau_{\mathcal{O}_{\mu_1}})^\ast\cdot\omega^L_{2\mathcal{O}_{\mu_2}}\\
& =(\tau_{\mathcal{O}_{\mu_2}}\cdot (T\varphi)_{\mathcal{O}_\mu})^\ast
\cdot\omega^L_{2\mathcal{O}_{\mu_2}}\\
&=((T\varphi)_{\mathcal{O}_\mu})^\ast\cdot
\tau_{\mathcal{O}_{\mu_2}}^\ast
\cdot\omega^L_{2\mathcal{O}_{\mu_2}}\\ &
=(j_{\mathcal{O}_{\mu_2}}^{-1}\cdot T\varphi \cdot
j_{\mathcal{O}_{\mu_1}})^\ast\cdot
j_{\mathcal{O}_{\mu_2}}^\ast\cdot\omega^L_2-(T\varphi)_{\mathcal{O}_\mu}^\ast
\cdot (\mathbf{J}_L)_{2\mathcal{O}_{\mu_2}}^\ast \cdot\omega^{L+}_{2\mathcal{O}_{\mu_2}}\\
&=j_{\mathcal{O}_{\mu_1}}^\ast\cdot(T\varphi)^\ast\cdot\omega^L_2
-(\mathbf{J}_L)_{1\mathcal{O}_{\mu_1}}^\ast\cdot\omega^{L+}_{1\mathcal{O}_{\mu_1}}\\
& =j_{\mathcal{O}_{\mu_1}}^\ast\cdot\omega^L_1
-(\mathbf{J}_L)_{1\mathcal{O}_{\mu_1}}^\ast\cdot\omega^{L+}_{1\mathcal{O}_{\mu_1}}\\
& = \tau_{\mathcal{O}_{\mu_1}}^\ast\cdot\omega^L_{1\mathcal{O}_{\mu_1}}.
\end{align*}
Because $\tau_{\mathcal{O}_{\mu_1}}$ is surjective, thus,
$((T\varphi)_{\mathcal{O}_\mu/G})^\ast\cdot\omega^L_{2\mathcal{O}_{\mu_2}}
=\omega^L_{1\mathcal{O}_{\mu_1}}$. \\

Conversely, assume that the $R_o$-reduced Lagrangian systems $((T
Q_i)_{\mathcal{O}_{\mu_i}}, \omega^L_{i\mathcal{O}_{\mu_i}},
l_{i\mathcal{O}_{\mu_i}})$, $i=1,2,$ are equivalent,
then there exists a diffeomorphism
$(T\varphi)_{\mathcal{O}_\mu/G}:(TQ_1)_{\mathcal{O}_{\mu_1}}\rightarrow
(TQ_2)_{\mathcal{O}_{\mu_2}}$, which is symplectic with respect to
the $R_o$-reduced Lagrangian symplectic forms $\omega^L_{i\mathcal{O}_i}, \;
i=1,2$, that is, $(T\varphi)_{\mathcal{O}_\mu/G}^\ast\cdot
\omega^L_{2\mathcal{O}_{\mu_2}}=\omega^L_{1\mathcal{O}_{\mu_1}}$.
Thus, from the above commutative Diagram-6, we can define a map
$(T\varphi)_{\mathcal{O}_\mu}:
(\mathbf{J}_L)_1^{-1}(\mathcal{O}_{\mu_1})\rightarrow
(\mathbf{J}_L)_2^{-1}(\mathcal{O}_{\mu_2}), $ such that
$\tau_{\mathcal{O}_{\mu_2}}\cdot
(T\varphi)_{\mathcal{O}_\mu}=(T\varphi)_{\mathcal{O}_\mu/G} \cdot
\tau_{\mathcal{O}_{\mu_1}}, $ and map $T\varphi: TQ_1\rightarrow
TQ_2, $ such that $j_{\mathcal{O}_{\mu_2}}\cdot
(T\varphi)_{\mathcal{O}_\mu}= T\varphi\cdot j_{\mathcal{O}_{\mu_1}},
$, as well as a diffeomorphism
$\varphi: Q_1\rightarrow Q_2$, whose tangent map is just $T\varphi:
TQ_1\rightarrow TQ_2 $. From definition of
$(T\varphi)_{\mathcal{O}_\mu}$ we know that
$(T\varphi)_{\mathcal{O}_\mu}$ is $(G_1,G_2)$-equivariant.\\

Now we shall show that $T\varphi$ is symplectic with respect to
the Lagrangian symplectic forms $\omega^L_i, \; i=1,2$, that is,
$\omega^L_1=(T\varphi)^\ast\cdot\omega^L_2$. In fact, since
$(T\varphi)_{\mathcal{O}_\mu/G}:
(TQ_1)_{\mathcal{O}_{\mu_1}}\rightarrow
(TQ_2)_{\mathcal{O}_{\mu_2}}$ is symplectic with respect to their
$R_o$-reduced Lagrangian symplectic forms, the map
$((T\varphi)_{\mathcal{O}_\mu/G})^\ast:
\Omega^2((TQ_2)_{\mathcal{O}_{\mu_2}})\rightarrow \Omega^2((T
Q_1)_{\mathcal{O}_{\mu_1}})$ satisfies
$((T\varphi)_{\mathcal{O}_\mu/G})^\ast\cdot
\omega^L_{2\mathcal{O}_{\mu_2}}=\omega^L_{1\mathcal{O}_{\mu_1}}$.
From (2.4) we have that
$j_{\mathcal{O}_{\mu_i}}^\ast\cdot\omega^L_i=\tau_{\mathcal{O}_{\mu_i}}^\ast\cdot
\omega^L_{i\mathcal{O}_{\mu_i}}
+(\mathbf{J}_L)_{i\mathcal{O}_{\mu_i}}^\ast\cdot\omega^{L+}_{i\mathcal{O}_{\mu_i}}$,
$i=1,2$, from the commutative Diagram-7, we have that
\begin{align*}
j_{\mathcal{O}_{\mu_1}}^\ast\cdot\omega^L_1 &
=\tau_{\mathcal{O}_{\mu_1}}^\ast\cdot
\omega^L_{1\mathcal{O}_{\mu_1}}+
(\mathbf{J}_L)_{1\mathcal{O}_{\mu_1}}^\ast\cdot\omega^{L+}_{1\mathcal{O}_{\mu_1}}\\
& =\tau_{1\mathcal{O}_{\mu_1}}^\ast\cdot
((T\varphi)_{\mathcal{O}_\mu/G})^\ast\cdot\omega^L_{2\mathcal{O}_{\mu_2}}
+(\mathbf{J}_L)_{1\mathcal{O}_{\mu_1}}^\ast\cdot\omega^{L+}_{1\mathcal{O}_{\mu_1}}\\
&=((T\varphi)_{\mathcal{O}_\mu/G}\cdot
\tau_{\mathcal{O}_{\mu_1}})^\ast\cdot
\omega^L_{2\mathcal{O}_{\mu_2}}
+(\mathbf{J}_L)_{1\mathcal{O}_{\mu_1}}^\ast\cdot\omega^{L+}_{1\mathcal{O}_{\mu_1}}\\
& =(\tau_{\mathcal{O}_{\mu_2}}\cdot
(T\varphi)_{\mathcal{O}_\mu})^\ast\cdot
\omega^L_{2\mathcal{O}_{\mu_2}}+ (\mathbf{J}_L)_{1\mathcal{O}_{\mu_1}}^\ast\omega^{L+}_{1\mathcal{O}_{\mu_1}}\\
&=(j_{\mathcal{O}_{\mu_2}}^{-1}\cdot T\varphi \cdot
j_{\mathcal{O}_{\mu_1}})^\ast \cdot
\tau_{\mathcal{O}_{\mu_2}}^\ast\cdot\omega^L_{2\mathcal{O}_{\mu_2}}
+(\mathbf{J}_L)_{1\mathcal{O}_{\mu_1}}^\ast\cdot\omega^{L+}_{1\mathcal{O}_{\mu_1}}\\
&=j_{\mathcal{O}_{\mu_1}}^\ast\cdot(T\varphi)^\ast\cdot
(j_{\mathcal{O}_{\mu_2}}^{-1})^\ast \cdot
[j_{\mathcal{O}_{\mu_2}}^{\ast}\cdot\omega^L_2
-(\mathbf{J}_L)_{2\mathcal{O}_{\mu_2}}^{\ast}\cdot\omega^{L+}_{2\mathcal{O}_{\mu_2}}]
+(\mathbf{J}_L)_{1\mathcal{O}_{\mu_1}}^{\ast}\cdot\omega^{L+}_{1\mathcal{O}_{\mu_1}}\\
&=j_{\mathcal{O}_{\mu_1}}^{\ast}\cdot(T\varphi)^\ast\cdot\omega^L_2-
((T\varphi)_{\mathcal{O}_\mu})^\ast \cdot
(\mathbf{J}_L)_{2\mathcal{O}_{\mu_2}}^{\ast}\cdot\omega^{L+}_{2\mathcal{O}_{\mu_2}}
+(\mathbf{J}_L)_{1\mathcal{O}_{\mu_1}}^{\ast}\cdot\omega^{L+}_{1\mathcal{O}_{\mu_1}}
\end{align*}
Notice that $j_{\mathcal{O}_{\mu_1}}$ is injective, and by
our hypothesis,
$$(\mathbf{J}_L)_{1\mathcal{O}_{\mu_1}}^{\ast}\cdot\omega^{L+}_{1\mathcal{O}_{\mu_1}}
=((T\varphi)_{\mathcal{O}_\mu})^\ast \cdot
(\mathbf{J}_L)_{2\mathcal{O}_{\mu_2}}^{\ast}\cdot\omega^{L+}_{2\mathcal{O}_{\mu_2}},$$
then $\omega^L_1=(T\varphi)^\ast \omega^L_2$.
Thus, the regular orbit reducible Lagrangian systems $(TQ_i, G_i,
\omega^L_i, L_i )$, $i=1,2,$ are equivalent.  \hskip 1cm
$\blacksquare $ \\

Thus, the regular orbit reduction Theorem 5.3 for the RCL systems
can be regarded as an extension of the regular
orbit reduction Theorem 5.4 for the regular Lagrangian systems under regular
controlled Lagrangian equivalence conditions.

\begin{rema}
If $(TQ, \omega^L)$ is a connected symplectic manifold, and
$\mathbf{J}_L: TQ\rightarrow \mathfrak{g}^\ast$ is a non-equivariant
momentum map with a non-equivariance group one-cocycle $\sigma:
G\rightarrow \mathfrak{g}^\ast$, which is defined by
$\sigma(g):=\mathbf{J}_L(g\cdot
z)-\operatorname{Ad}^\ast_{g^{-1}}\mathbf{J}_L(z)$, where $g\in G$
and $z\in TQ$. Then we know that $\sigma$ produces a new affine
action $\Theta: G\times \mathfrak{g}^\ast \rightarrow
\mathfrak{g}^\ast $ defined by
$\Theta(g,\mu):=\operatorname{Ad}^\ast_{g^{-1}}\mu + \sigma(g)$,
where $\mu \in \mathfrak{g}^\ast$, with respect to which the given
momentum map $\mathbf{J}_L$ is equivariant. Assume that $G$ acts
freely and properly on $TQ$, and $\mathcal{O}_\mu= G\cdot \mu
\subset \mathfrak{g}^\ast$ denotes the G-orbit of the point $\mu \in
\mathfrak{g}^\ast$ with respect to the above affine action $\Theta$,
and $\mu$ is a regular value of $\mathbf{J}_L$. Then the quotient
space $(TQ)_{\mathcal{O}_\mu}=\mathbf{J}_L^{-1}(\mathcal{O}_\mu)/ G
$ is also a symplectic manifold with the symplectic form
$\omega^L_{\mathcal{O}_\mu}$ uniquely characterized by $(2.4)$. In
this case, we can also define the regular orbit reducible RCL system
$(TQ,G,\omega^L,L,F^L,\mathcal{C}^L)$ and RoCL-equivalence, and
prove the regular orbit reduction theorem for the RCL system by using
the above similar way.
\end{rema}

It is worthy of noting that the research idea and work in Marsden et al. 
\cite{mawazh10} are very important. 
The authors not only correct and renew carefully some wrong definitions
for  CH system and its reduced CH systems, as well as
CH-equivalence and the reduced CH-equivalence
in Chang et al. \cite{chbllemawo02, chma04},
but also set up the regular reduction theory of regular
controlled Hamiltonian systems on a symplectic fiber bundle, 
from the viewpoint of completeness of Marsden-Weinstein reduction.
In this paper, following the ideas in Marsden et al. \cite{mawazh10},
we correct and renew carefully some wrong definitions
for  CL system and its reduced CL systems, as well as
CL-equivalence and the reduced CL-equivalence
in Chang et al. \cite{chbllemawo02, chma04},
and set up the regular reduction theory of regular
controlled Lagrangian systems on a symplectic fiber bundle,
by analyzing carefully the geometrical and topological
structures of the phase space and the reduced phase space of the
regular Lagrangian system. Note that 
some developments around the work in Marsden et al. \cite{mawazh10}
are given in Wang and Zhang \cite{wazh12}, Ratiu and Wang \cite{rawa12}, 
Wang \cite{wa18}, and Wang \cite{wa17}, 
and some applications are given in Wang \cite{wa20a, wa13e}. 
Thus, it is natural idea to develop a variety of reduction theory
and applications for regular controlled Lagrangian systems,
in particular, in celestial mechanics, hydrodynamics and
plasma physics. In addition, it is also
an important topic for us to explore and reveal the deeply internal
relationships between the geometrical structures of phase spaces and the dynamical
vector fields of the controlled mechanical systems.
In particular, 
it is an important task for us to correct and develop well
the research work of Professor Jerrold E. Marsden,
such that we never feel sorry for his great cause.\\

\end{document}